\documentclass[11pt]{article}
\usepackage[top=1in, bottom=1in, left=1in, right=1in]{geometry}
\usepackage{tikz}
\usetikzlibrary{shapes,arrows,matrix,patterns,positioning,calc}
\usepackage{color}
\usepackage{graphicx}
\usepackage{algorithmic}
\usepackage{amssymb}
\usepackage{epstopdf}
\usepackage[ruled]{algorithm2e}
\usepackage{float}
\usepackage{amsmath,amsthm,bm,color,epsfig,enumerate,caption}
\usepackage{hyperref}
\usepackage{booktabs}
\usepackage{multirow}
\usepackage{array}

\newtheorem{theorem}{Theorem}[section]

\renewcommand{\appendix}[1]{\section*{Appendix: #1}}
\newcommand{\norm}[1]{\left\lVert#1\right\rVert}
\DeclareMathOperator{\dist}{dist}

\renewcommand{\O}{O}
\renewcommand{\i}{i}
\renewcommand{\j}{j}
\newcommand{\I}{\mathcal{I}}

\newcommand{\bbC}{\mathbb{C}}
\newcommand{\bbR}{\mathbb{R}}
\newcommand{\ve}[1]{{#1}}

\newcommand{\vcenteredinclude}[1]{\begingroup
\setbox0=\hbox{\includegraphics[width=3.5cm]{#1}}
\parbox{\wd0}{\box0}\endgroup}

\makeatletter
\newcommand*{\extendadd}{
  \mathbin{
    \mathpalette\extend@add{}
  }
}
\newcommand*{\extend@add}[2]{
  \ooalign{
    $\m@th#1\leftrightarrow$
    \vphantom{$\m@th#1\updownarrow$}
    \cr
    \hfil$\m@th#1\updownarrow$\hfil
  }
}
\makeatother

\begin{document}
\title{Multidimensional Butterfly Factorization}

\author{Yingzhou Li$^\sharp$,
  Haizhao Yang$^*$,
  Lexing Ying$^{\dagger\sharp}$
  \vspace{0.1in}\\
  $\dagger$ Department of Mathematics, Stanford University\\
  $\sharp$ ICME, Stanford University\\
  $*$ Department of Mathematics, Duke University
}

\maketitle

\begin{abstract}
This paper introduces the multidimensional butterfly factorization as
a data-sparse representation of multidimensional kernel matrices that
satisfy the complementary low-rank property. This factorization
approximates such a kernel matrix of size $N\times N$ with a product
of $\O(\log N)$ sparse matrices, each of which contains $\O(N)$
nonzero entries.  We also propose efficient algorithms for
constructing this factorization when either (i) a fast algorithm for
applying the kernel matrix and its adjoint is available or (ii) every
entry of the kernel matrix can be evaluated in $\O(1)$ operations. For
the kernel matrices of multidimensional Fourier integral operators,
for which the complementary low-rank property is not satisfied due to
a singularity at the origin, we extend this factorization by combining
it with either a polar coordinate transformation or a multiscale
decomposition of the integration domain to overcome the singularity.
Numerical results are provided to demonstrate the efficiency of the
proposed algorithms.
\end{abstract}




\section{Introduction}
\label{sec:intro}

\subsection{Problem statement}

This paper is concerned with the efficient evaluation of
\begin{equation}\label{eq:kernel}
  u(x) = \sum_{\xi\in \Omega} K(x,\xi)g(\xi),\quad x\in X,
\end{equation}
where $X$ and $\Omega$ are typically point sets in $\bbR^d$ for $d\geq
2$, $K(x,\xi)$ is a kernel function that satisfies a complementary
low-rank property, $g(\xi)$ is an input function for $\xi\in \Omega$,
and $u(x)$ is an output function for $x\in X$.  To define this
complementary low-rank property for multidimensional kernel matrices,
we first assume that without loss of generality there are $N$ points
in each point set. In addition, the domains $X$ and $\Omega$ are
associated with two hierarchical trees $T_X$ and $T_\Omega$,
respectively, where each node of these trees represents a subdomain of
$X$ or $\Omega$. Both $T_X$ and $T_\Omega$ are assumed to have
$L=\O(\log N)$ levels with $X$ and $\Omega$ being the roots at level
$0$. The computation of \eqref{eq:kernel} is essentially a matrix
vector multiplication
\[
u = K g,
\]
where $K := (K(x,\xi))_{x\in X,\xi\in\Omega}$,
$g:=(g(\xi))_{\xi\in\Omega}$, and $u := (u(x))_{x\in X}$ by a slight
abuse of notations. The matrix $K$ is said to satisfy the {\em
  complementary low-rank property} if for any level $\ell$ between $0$
and $L$ and for any node $A$ on the $\ell$-th level of $T_X$ and any
node $B$ on the $(L-\ell)$-th level of $T_\Omega$, the submatrix
$K_{A,B}:=(K(x_i,\xi_j))_{x_i\in A, \xi_j\in B}$ is numerically
low-rank with the rank bounded by a uniform constant independent of
$N$.  In most applications, this numerical rank is bounded
polynomially in $\log(1/\epsilon)$ for a given precision $\epsilon$. A
well-known example of such a matrix is the multidimensional Fourier
transform matrix.

For a complementary low-rank kernel matrix $K$, the {\em butterfly
  algorithm} developed in \cite{fio07,fio09,mmd,1dba,wavemoth} enables
one to evaluate the matrix-vector multiplication in $O(N\log N)$
operations. More recently in \cite{1dbf}, we introduced the {\em
  butterfly factorization} as a data-sparse multiplicative
factorization of the kernel matrix $K$ in the one-dimensional case
($d=1$):
\begin{equation}\label{eq:GBF}
  K\approx U^LG^{L-1}\cdots G^{L/2}M^{L/2}
  \left(H^{L/2}\right)^*\cdots\left(H^{L-1}\right)^*\left(V^L\right)^*,
\end{equation}
where the depth $L=\O(\log N)$ is assumed to be an even number and
every factor in \eqref{eq:GBF} is a sparse matrix with $\O(N)$
nonzero entries. Here the superscript of a matrix denotes the level of
the factor rather than the power of a matrix. This factorization
requires $\O(N\log N)$ memory and applying \eqref{eq:GBF} to any
vector takes $\O(N\log N)$ operations once the factorization is
computed. In fact, one can view the factorization in \eqref{eq:GBF} as
a compact algebraic representation of the butterfly algorithm. In
\cite{1dbf}, we also introduced algorithms for constructing the
butterfly factorization for the following two cases:
\begin{enumerate}[(i)]
\item A black-box routine for rapidly computing $Kg$ and $K^*g$ in
    $\O(N\log N)$ operations is available;
\item A routine for evaluating any entry of $K$ in $\O(1)$ operations
  is given.
\end{enumerate}
In this paper, we turn to the butterfly factorization for the
multidimensional problems and describe how to construct them for these
two cases.

When the kernel strictly satisfies the complementary low-rank property
(e.g., the non-uniform FFT), the algorithms proposed in \cite{1dbf} can
be generalized in a rather straightforward way. This is presented in
detail in Section \ref{sec:gbf}.

However, many important multidimensional kernel matrices fail to
satisfy the complementary low-rank property in the entire domain
$X\times\Omega$. Among them, the most significant example is probably
the Fourier integral operator, which typically has a singularity at
the origin $\xi=0$ in the $\Omega$ domain. For such an example,
existing butterfly algorithms provide two solutions.
\begin{itemize}
\item The first one, proposed in \cite{fio09}, removes the singularity
  by applying a polar transformation that maps the domain $\Omega$
  into a new domain $P$. After this transformation, the new kernel
  matrix defined on $X\times P$ satisfies the complementary low-rank
  property and one can then apply the butterfly factorization in the
  $X$ and $P$ domain instead. This is discussed in detail in Section
  \ref{sec:pbf} and we refer to this algorithm as the {\em polar
    butterfly factorization} (PBF).
\item The second solution proposed in \cite{mba} is based on the
  observation that, though not on the entire $\Omega$ domain, the
  complementary low-rank property holds in subdomains of $\Omega$ that
  are well separated from the origin in a certain sense.  For example,
  one can start by partitioning the domain $\Omega$ into a disjoint
  union of a small square $\Omega_C$ covering $\xi=0$ and a sequence
  of dyadic coronas $\Omega_t$, i.e., $\Omega = \Omega_C \cup
  \left(\cup_t\Omega_t\right)$. Accordingly, one can rewrite the
  kernel evaluation \eqref{eq:kernel} as a summation of the form
  \begin{equation}
    \label{eq:Ksum}
    K = K_C R_C +\sum_t K_t R_t,
  \end{equation}
  where $K_C$ and $K_t$ are the kernel matrices restricted to $X\times
  \Omega_C$ and $X\times \Omega_t$, $R_C$ and $R_t$ are the operators
  of restricting the input functions defined on $\Omega$ to the
  subdomain $\Omega_C$ and $\Omega_t$, respectively. In fact, each
  kernel $K_t$ satisfies the complementary low-rank property and hence
  one can approximate it with the multidimensional butterfly
  factorization in Section \ref{sec:gbf}. Combining the factorizations
  for all $K_t$ with \eqref{eq:Ksum} results the {\em multiscale
    butterfly factorization} (MBF) for the entire matrix $K$ and this
  will be discussed in detail in Section \ref{sec:mbf}.
\end{itemize}

In order to simplify the presentation, this paper focuses on the two
dimensional case ($d=2$). Furthermore, we assume that the points in $X$
and $\Omega$ are uniformly distributed in both domains as follows:
\begin{equation}
  \label{eq:X}
  X = \left\{ x = \left( \frac{n_1}{n}, \frac{n_2}{n}\right),
      0 \leq  n_1,n_2 < n\text{ with }
  n_1, n_2 \in \mathbb{Z} \right\}
\end{equation}
and
\begin{equation}
  \label{eq:Omega}
  \Omega = \left\{ \xi = (n_1, n_2),
      - \frac{n}{2} \leq n_1,n_2 < \frac{n}{2}\text{ with }
  n_1, n_2 \in \mathbb{Z} \right\},
\end{equation}
where $n$ is the number of points in each dimension and $N=n^2$.
This is the standard setup for two dimensional Fourier transforms and
FIOs. 

\subsection{Related work}
\label{sec:moti}

For a complementary low-rank kernel matrix $K$, the butterfly
algorithm provides an efficient way for evaluating
\eqref{eq:kernel}. It was initially proposed in \cite{mmd} and further
developed in \cite{fio09,hu,mba,1dba,fio13,wavemoth,sht,sft}. One can
roughly classify the existing butterfly algorithms into two
groups.
\begin{itemize}
\item The first group (e.g. \cite{1dba,wavemoth,sht}) requires a
  precomputation stage for constructing the low-rank
  approximations of the numerically low-rank submatrices of
  \eqref{eq:kernel}. This precomputation stage typically takes
  $\O(N^2)$ operations and uses $\O(N\log N)$ memory. Once the
  precomputation is done, the evaluation of \eqref{eq:kernel} can be
  carried out in $\O(N\log N)$ operations.
\item The second group (e.g. \cite{fio09,hu,mba,fio13}) assumes prior
  knowledge of analytic properties of the kernel function.  Under such
  analytic assumptions, one avoids precomputation by writing down the
  low-rank approximations for the numerically low-rank submatrices
  explicitly.  These algorithms typically evaluate \eqref{eq:kernel}
  with $\O(N\log N)$ operations.
\end{itemize}
In a certain sense, the algorithms proposed in this paper can be
viewed as a compromise of these two types. On the one hand, it makes
rather weak assumptions about the kernel. Instead of requiring the
kernel function as was done for the second type, we only assume that
either (i) a fast matrix-vector multiplication routine or (ii) a
kernel matrix sampling routine is available. On the other hand, these
new algorithms reduce the precomputation cost to $O(N^{3/2} \log N)$,
as compared to the quadratic complexity of the first group.

The multidimensional butterfly factorization can also be viewed as a
process of recovering a structured matrix via either sampling or
matrix-vector multiplication. There has been a sequence of articles in
this line of research. For example, we refer to
\cite{randsvd,Rec2,Rec3} for recovering numerically low-rank matrices,
\cite{HSSMatrix} for recovering an {HSS}
matrices, and \cite{HMatrix} for recovering $\mathcal{H}$-matrices. This
paper generalizes the work of \cite{1dbf} by considering complementary
low-rank matrices coming from multidimensional problems.

\subsection{Organization}
\label{sec:orga}

The rest of this paper is organized as follows. Section \ref{sec:gbf}
reviews the basic tools and describes the {\em multidimensional
  butterfly factorization} for kernel matrices that strictly satisfy
the complementary low-rank property. We then extend it in two
different ways to address the multidimensional Fourier integral
operators. Section \ref{sec:pbf} introduces the polar butterfly
factorization (PBF) based on the polar butterfly algorithm proposed in
\cite{fio09}. Section \ref{sec:mbf} discusses the multiscale butterfly
factorization (MBF) based on the multiscale butterfly algorithm
proposed in \cite{mba}. Finally, in Section \ref{sec:conc},
we conclude with some discussions.


\section{Two-Dimensional {Butterfly Factorization}}
\label{sec:gbf}

This section presents the two-dimensional butterfly factorization for
a kernel matrix $K = (K(x,\xi))_{x\in X,\xi\in \Omega}$ that satisfies
the complementary low-rank property in $X\times \Omega$ with $X$ and
$\Omega$ given in \eqref{eq:X} and \eqref{eq:Omega}.

\subsection{Randomized low-rank factorization}
\label{sec:randlr}

The butterfly factorization relies heavily on randomized procedures
for computing low-rank factorizations. For a matrix $Z\in
\bbC^{m\times n}$, a rank-$r$ approximation in 2-norm can be computed
via the truncated singular value decomposition ({SVD}),
\begin{equation}
  \label{eq:Z}
  Z\approx U_0\Sigma_0V_0^*,
\end{equation}
where $U_0\in \bbC^{m\times r}$ and $V_0\in \bbC^{n\times r}$ are
unitary matrices, $\Sigma_0\in \bbR^{r\times r}$ is a diagonal matrix
with the largest $r$ singular values of $Z$ in decreasing order.

Once $Z\approx U_0\Sigma_0V_0^*$ is available, we can also construct
different low-rank factorizations of $Z$ in three forms:
\begin{align}
  & Z\approx USV^*,\quad U=U_0\Sigma_0,\quad S = \Sigma_0^{-1},
  \quad V^*=\Sigma_0V_0^*;\label{eq:lowrank1}\\
  & Z\approx UV^*,\quad U=U_0\Sigma_0, \quad V^*=V_0^*;\label{eq:lowrank2}\\
  & Z\approx UV^*,\quad U=U_0,  \quad V^*=\Sigma_0V_0^*.\label{eq:lowrank3}
\end{align}
As we shall see, the butterfly factorization uses each of these three
forms in different stages of the algorithm.

In \cite{1dbf}, we showed that the rank-$r$ {SVD} \eqref{eq:Z} can be
constructed approximately via either random matrix-vector
multiplication \cite{randsvd} or random sampling
\cite{randsamp1,randsamp2}. In both cases, the key is to find accurate
approximate bases for both the column and row spaces of $Z$ and
approximate the largest $r$ singular values using these bases.

\paragraph{SVD via random matrix-vector multiplication.} 
This algorithm proceeds as follows.
\begin{itemize}
\item This algorithm first applies $Z$ to a Gaussian random matrix
  $C\in \bbC^{n\times (r+k)}$ and its adjoint $Z^*$ to a Gaussian
  random matrix $R\in \bbC^{m\times (r+k)}$, where $k$ is the
  oversampling constant.
\item Second, computing the pivoted {QR} decompositions of $ZC$ and
  $Z^*R$ identifies unitary matrices $Q_{col}\in \bbC^{m\times r}$ and
  $Q_{row}\in \bbC^{n\times r}$, which approximately span the column
  and row spaces of $Z$, respectively. 
\item Next, the algorithms seeks a matrix $M$ that satisfies
  \[
  Z \approx Q_{col} M Q_{row}^*
  \] 
  by setting $M = (R^*Q_{col})^\dagger R^*ZC(Q_{row}^*C)^\dagger$,
  where $(\cdot)^\dagger$ denotes the pseudo inverse.  
\item Finally, combining the singular value decomposition $M=U_M
  \Sigma_M V_M^*$ of the matrix $M$ with the above approximation
  results in the desired approximate rank-$r$ SVD
  \[
  Z \approx (Q_{col} U_M) \Sigma_M (Q_{row} V_M)^*.
  \]
\end{itemize}
Suppose that the cost of applying $Z$ and $Z^*$ to an arbitrary vector
is $C_Z(m,n)$. Then the construction complexity of this procedure is
$\O(C_Z(m,n) r+\max(m,n)r^2)$. As we shall see later, when the
black-box routines for rapidly applying $K$ and $K^*$ are available,
this procedure would be embedded into the algorithms for constructing
the butterfly factorizations.

\paragraph{SVD via random sampling.} This algorithm proceeds as follows.
\begin{itemize}
\item The first stage discovers the representative columns and rows
  progressively via computing multiple pivoted {QR} factorizations on
  randomly selected rows and columns of $Z$.  The representative
  columns and rows are set to be empty initially. As the procedure
  processes, more and more columns (rows) are marked as {\em
    representative} and they are used in turn to discover new
  representative rows (columns). The procedure stops when the sets of
  the representative rows and columns stabilize. At this point, the
  representative columns (rows) approximately span the column (row)
  spaces of $Z$.
\item Second, computing the pivoted {QR} decompositions of the
  representative columns and rows identifies unitary matrices
  $Q_{col}\in \bbC^{m\times r}$ and $Q_{row}\in \bbC^{n\times r}$,
  which approximately span the column and row spaces of $Z$,
  respectively.
\item Next, the algorithm seeks a matrix $M$ that satisfies
  \[
  Z \approx Q_{col} M Q_{row}^*.
  \] 
  This is done by restricting this equation to a random row set
  $I_{row}$ and a random column set $I_{col}$ and consider
  \[
  Z(I_{row},I_{col}) \approx Q_{col}(I_{row},:) M Q_{row}(I_{col},:)^*.
  \]
  Here both $I_{row}$ and $I_{col}$ are of size $O(r)$ and we require
  $I_{row}$ and $I_{col}$ to contain the set of representative rows and
  columns, respectively. From the above equation, we can solve $M$ by
  setting
  \[
  M = (Q_{col}(I_{row},:))^\dagger Z(I_{row},I_{col}) (Q_{row}(I_{col},:)^*)^\dagger.
  \]
\item Finally, combining the singular value decomposition $M=U_M
  \Sigma_M V_M^*$ of the matrix $M$ with the approximation $Z \approx
  Q_{col} M Q_{row}^*$ results in the desired approximate rank-$r$
  SVD
  \[
  Z \approx (Q_{col} U_M) \Sigma_M (Q_{row} V_M)^*.
  \]
\end{itemize}
The construction complexity of this procedure is $\O(\max(m,n)r^2)$
in practice. When an arbitrary entry of $Z$ can be evaluated in
$\O(1)$ operations, this procedure is the method of choice for
constructing low-rank factorizations.

\subsection{Notations and overall structure}

We adopt the notation of the one-dimensional butterfly factorization
introduced in \cite{1dbf} and adjust them to the two-dimensional case
of this paper.

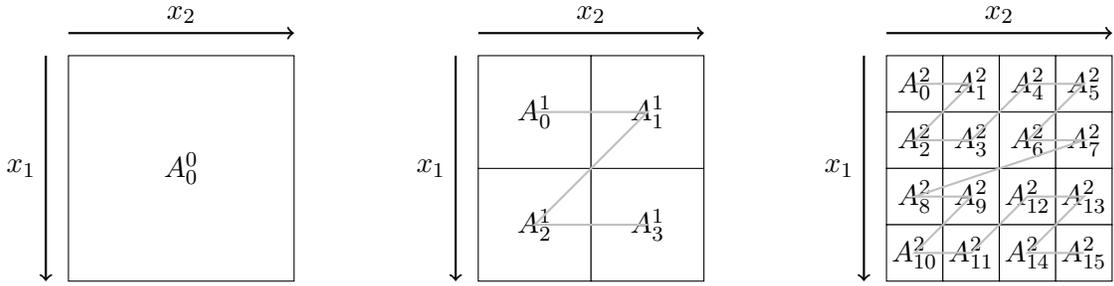
\begin{figure}[htb]
\begin{minipage}{.33\textwidth}
\centering
\begin{tikzpicture}[scale=3]
\fill[white] (0,0) rectangle (1,1);
\draw[black] (0,0) rectangle (1,1);
\draw (0.5,0.5) node[rectangle] {$A^0_{{0}}$};
\draw[->,thick] (0,1.1) -> (1,1.1);
\draw (0.5,1.1) node[rectangle,above] {$x_2$};
\draw[->,thick] (-0.1,1) -> (-0.1,0);
\draw (-0.1,0.5) node[rectangle,left] {$x_1$};
\end{tikzpicture}
\end{minipage}%
\begin{minipage}{.33\textwidth}
\centering
\begin{tikzpicture}[scale=3]
\fill[white] (0,0) rectangle (0.5,0.5);
\draw[black] (0,0) rectangle (0.5,0.5);
\draw (0.25,0.25) node[rectangle] {$A^1_{2}$};
\fill[white] (0,0.5) rectangle (0.5,1);
\draw[black] (0,0.5) rectangle (0.5,1);
\draw (0.25,0.75) node[rectangle] {$A^1_{0}$};
\fill[white] (0.5,0) rectangle (1,0.5);
\draw[black] (0.5,0) rectangle (1,0.5);
\draw (0.75,0.25) node[rectangle] {$A^1_{3}$};
\fill[white] (0.5,0.5) rectangle (1,1);
\draw[black] (0.5,0.5) rectangle (1,1);
\draw (0.75,0.75) node[rectangle] {$A^1_{1}$};
\draw[lightgray,thick] (0.25,0.75) -- (0.75,0.75) -- (0.25,0.25) -- (0.75,0.25);
\draw[->,thick] (0,1.1) -> (1,1.1);
\draw (0.5,1.1) node[rectangle,above] {$x_2$};
\draw[->,thick] (-0.1,1) -> (-0.1,0);
\draw (-0.1,0.5) node[rectangle,left] {$x_1$};
\end{tikzpicture}
\end{minipage}%
\begin{minipage}{.33\textwidth}
\centering
\begin{tikzpicture}[scale=3]
\fill[white] (0,0.75) rectangle (0.25,1);
\draw[black] (0,0.75) rectangle (0.25,1);
\draw (0.125,0.875) node[rectangle] {$A^2_{0}$};
\fill[white] (0.25,0.75) rectangle (0.5,1);
\draw[black] (0.25,0.75) rectangle (0.5,1);
\draw (0.375,0.875) node[rectangle] {$A^2_{1}$};
\fill[white] (0,0.5) rectangle (0.25,0.75);
\draw[black] (0,0.5) rectangle (0.25,0.75);
\draw (0.125,0.625) node[rectangle] {$A^2_{2}$};
\fill[white] (0.25,0.5) rectangle (0.5,0.75);
\draw[black] (0.25,0.5) rectangle (0.5,0.75);
\draw (0.375,0.625) node[rectangle] {$A^2_{3}$};
\fill[white] (0.5,0.75) rectangle (0.75,1);
\draw[black] (0.5,0.75) rectangle (0.75,1);
\draw (0.625,0.875) node[rectangle] {$A^2_{4}$};
\fill[white] (0.75,0.75) rectangle (1,1);
\draw[black] (0.75,0.75) rectangle (1,1);
\draw (0.875,0.875) node[rectangle] {$A^2_{5}$};
\fill[white] (0.5,0.5) rectangle (0.75,0.75);
\draw[black] (0.5,0.5) rectangle (0.75,0.75);
\draw (0.625,0.625) node[rectangle] {$A^2_{6}$};
\fill[white] (0.75,0.5) rectangle (1,0.75);
\draw[black] (0.75,0.5) rectangle (1,0.75);
\draw (0.875,0.625) node[rectangle] {$A^2_{7}$};
\fill[white] (0,0.25) rectangle (0.25,0.5);
\draw[black] (0,0.25) rectangle (0.25,0.5);
\draw (0.125,0.375) node[rectangle] {$A^2_{8}$};
\fill[white] (0.25,0.25) rectangle (0.5,0.5);
\draw[black] (0.25,0.25) rectangle (0.5,0.5);
\draw (0.375,0.375) node[rectangle] {$A^2_{9}$};
\fill[white] (0,0) rectangle (0.25,0.25);
\draw[black] (0,0) rectangle (0.25,0.25);
\draw (0.125,0.125) node[rectangle] {$A^2_{10}$};
\fill[white] (0.25,0) rectangle (0.5,0.25);
\draw[black] (0.25,0) rectangle (0.5,0.25);
\draw (0.375,0.125) node[rectangle] {$A^2_{11}$};
\fill[white] (0.5,0.25) rectangle (0.75,0.5);
\draw[black] (0.5,0.25) rectangle (0.75,0.5);
\draw (0.625,0.375) node[rectangle] {$A^2_{12}$};
\fill[white] (0.75,0.25) rectangle (1,0.5);
\draw[black] (0.75,0.25) rectangle (1,0.5);
\draw (0.875,0.375) node[rectangle] {$A^2_{13}$};
\fill[white] (0.5,0) rectangle (0.75,0.25);
\draw[black] (0.5,0) rectangle (0.75,0.25);
\draw (0.625,0.125) node[rectangle] {$A^2_{14}$};
\fill[white] (0.75,0) rectangle (1,0.25);
\draw[black] (0.75,0) rectangle (1,0.25);
\draw (0.875,0.125) node[rectangle] {$A^2_{15}$};
\draw[lightgray, thick]
                 (0.125,0.875) -- (0.375,0.875) -- (0.125,0.625) -- (0.375,0.625)
              -- (0.625,0.875) -- (0.875,0.875) -- (0.625,0.625) -- (0.875,0.625)
              -- (0.125,0.375) -- (0.375,0.375) -- (0.125,0.125) -- (0.375,0.125)
              -- (0.625,0.375) -- (0.875,0.375) -- (0.625,0.125) -- (0.875,0.125);

\draw[->,thick] (0,1.1) -> (1,1.1);
\draw (0.5,1.1) node[rectangle,above] {$x_2$};
\draw[->,thick] (-0.1,1) -> (-0.1,0);
\draw (-0.1,0.5) node[rectangle,left] {$x_1$};
\end{tikzpicture}
\end{minipage}%
\caption{An illustration of Z-order curve cross levels. The
  superscripts indicate the different levels while the subscripts
  indicate the index in the Z-ordering. The light gray lines show the
  ordering among the subdomains on the same level.  Left: The root at
  level $0$.  Middle: At level $1$, the domain $A^0_{0}$ is divided
  into $2\times 2$ subdomains $A^1_\i$ with $\i\in \I^1=\{0,1,2,3\}$.
  These $4$ subdomains are ordered according to the Z-ordering. Right:
  At level $2$, the domain $A^0_0$ is divided into $4\times 4$
  subdomains $A^2_\i$ with $\i\in \I^2=\{0,1,\ldots,15\}$. These $16$
  subdomains are ordered similarly.}
\label{fig:domain-order-2D}
\end{figure}

Recall that $n$ is the number of grid points on each dimension and
$N=n^2$ is the total number of points. Suppose that $T_X$ and
$T_\Omega$ are complete quadtrees with $L = \log n$ levels and,
without loss of generality, $L$ is an even integer.  For a fixed level
$\ell$ between $0$ and $L$, the quadtree $T_X$ has $4^\ell$ nodes at
level $\ell$. By defining $\I^\ell = \{0,1,\ldots,4^\ell-1\}$, we
denote these nodes by $A^\ell_\i$ with $\i\in \I^\ell$. These $4^\ell$
nodes at level $\ell$ are further ordered according to a Z-order curve
(or Morton order) as illustrated in Figure
\ref{fig:domain-order-2D}. Based on this Z-ordering, the node
$A^\ell_\i$ at level $\ell$ has four child nodes denoted by
$A^{\ell+1}_{4\i + t}$ with $t=0,\dots,3$. The nodes plotted in Figure
\ref{fig:domain-order-2D} for $\ell = 1$ (middle) and $\ell=2$ (right)
illustrate the relationship between the parent node and its child
nodes. Similarly, in the quadtree $T_\Omega$, the nodes at the
$L-\ell$ the are denoted as $B^{L-\ell}_\j$ for $\j\in\I^{L-\ell}$.

For any level $\ell$ between $0$ and $L$, the kernel matrix $K$ can be
partitioned into $O(N)$ submatrices $K_{A^\ell_\i,B^{L-\ell}_\j}
:=(K(x,\xi))_{x\in A^\ell_\i,\xi\in B^{L-\ell}_\j}$ for $\i\in\I^\ell$
and $\j\in\I^{L-\ell}$. For simplicity, we shall denote
$K_{A^\ell_\i,B^{L-\ell}_\j}$ as $K^{\ell}_{\i,\j}$, where the
superscript $\ell$ denotes the level in the quadtree $T_X$. Because of
the complementary low-rank property, every submatrix $K^\ell_{\i,\j}$
is numerically low-rank with the rank bounded by a uniform constant
$r$ independent of $N$.

The two-dimensional butterfly factorization consists of two stages.
The first stage computes the factorizations
\[
K^h_{\i,\j}\approx U^h_{\i,\j}S^h_{\i,\j}\left(V^h_{\j,\i}\right)^*
\]
for all $\i,\j\in\I^h$ at the middle level $h=L/2$, following the form
\eqref{eq:lowrank1}. These factorizations can then be assembled into
three sparse matrices $U^h$, $M^h$, and $V^h$ to give rise to a
factorization for $K$:
\begin{equation}
  K\approx U^h M^h\left(V^h\right)^*.
\end{equation}
This stage is referred to as the \emph{middle level factorization} and
is described in Section \ref{sec:mlf}. In the second stage, we
recursively factorize the left and right factors $U^h$ and $V^h$ to
obtain
\begin{equation*}
  U^h \approx U^LG^{L-1}\cdots G^h \quad\text{and}\quad
  \left(V^h\right)^* \approx \left(H^h\right)^*\cdots
  \left(H^{L-1}\right)^*\left(V^L\right)^*,
\end{equation*}
where the matrices on the right hand side in each formula are sparse
matrices with $O(N)$ nonzero entries. Once they are ready, we
assemble all factors together to produce a data-sparse approximate
factorization for $K$:
\begin{equation}\label{eq:gbf}
  K\approx U^LG^{L-1}\cdots G^hM^h\left(H^h\right)^*\cdots
  \left(H^{L-1}\right)^*\left(V^L\right)^*,
\end{equation}
This stage is referred to as the \emph{recursive factorization} and is
discussed in Section \ref{sec:rf}.

\subsection{Middle level factorization}
\label{sec:mlf}

Recall that we consider the construction of multidimensional butterfly
factorization for two cases:
\begin{enumerate}[(i)]
\item A black-box routine for rapidly computing $Kg$ and $K^*g$ in
    $\O(N\log N)$ operations is available;
\item A routine for evaluating any entry of $K$ in $\O(1)$ operations
  is given.
\end{enumerate}

In Case (i), we construct an approximate rank-$r$ {SVD} of each
$K^h_{\i,\j}\in \bbR^{n\times n}$ with $\i,\j\in \I^h$ using the SVD
via random matrix-vector multiplication (the first option in
Section \ref{sec:randlr}). This requires applying each $K^h_{\i,\j}$ to
a {Gaussian} random matrix $C_\j\in\bbC^{n\times(r+k)}$ and its
adjoint to a Gaussian random matrix $R_\i\in\bbC^{(r+k)\times
  n}$. Here $r$ is the desired numerical rank and $k$ is the
oversampling parameter. If a black box routine for applying the matrix
$K$ and its adjoint is available, this can be done in an efficient way
as follows. For each $\j\in\I^h$, one constructs a zero-padded random
matrix $C^P_\j\in\bbC^{N\times (r+k)}$ by padding zero to $C_\j$. From
the relationship
\begin{equation}
  KC^P_\j = K
  \begin{pmatrix}
    0\\
    C_\j\\
    0
  \end{pmatrix}
  =
  \begin{pmatrix}
    K^h_{{0},\j}C_\j\\
    \vdots\\
    K^h_{{4^h-1},\j}C_\j
  \end{pmatrix},
\end{equation}
it is clear that applying $K$ to the matrix $C^P_\j$ produces
$K^h_{\i,\j}C_\j$ for all $\i\in\I^h$.  Similarly, we construct
zero-padded random matrices $R^P_\i\in \bbC^{N\times (r+k)}$ by padding
zero to $R_\i$ and compute
\begin{equation}
  K^* R^P_\i  =
  K^*
  \begin{pmatrix}
    0\\
    R_\i\\
    0
  \end{pmatrix}
  =
  \begin{pmatrix}
    \left( K^h_{\i,{0}}\right)^* R_\i  \\
    \vdots \\
    \left( K^h_{\i,{4^h-1}}\right)^* R_\i
  \end{pmatrix}
\end{equation}
by using the black-box routine for applying the adjoint of
$K$. Finally, the approximated rank-$r$ {SVD} of $K^h_{\i,\j}$ for
each pair of $\i\in\I^h$ and $\j\in\I^h$ is computed from
$K^h_{\i,\j}C_\j$ and $\left(K^h_{\i,\j}\right)^*R_\i$.

In Case (ii), since an arbitrary entry of $K$ can be evaluated in
$\O(1)$ operations, the approximate rank-$r$ {SVD} of $K^h_{\i,\j}$ is
computed using the SVD via randomized sampling
\cite{randsamp1,randsamp2} (the second option in Section
\ref{sec:randlr}).

In both cases, once the approximate rank-$r$ SVD is ready, we
transform it into the form of \eqref{eq:lowrank1}:
\begin{equation}
  \label{eq:Kij}
  K^h_{\i,\j}\approx U^h_{\i,\j}S^h_{\i,\j}\left(V^h_{\j,\i}\right)^*.
\end{equation}
Here the columns of the left and right factors $U^h_{\i,\j}$ and
$V^h_{\j,\i}$ are scaled by the singular values of $K^h_{\i,\j}$ such
that $U^h_{\i,\j}$ and $V^h_{\j,\i}$ keep track of the importance of
the column and row bases for further factorizations.

\begin{figure}[htp]
\begin{minipage}{\textwidth}
\centering
\resizebox{3.5cm}{!}{
\begin{tikzpicture}[baseline=-0.5ex]
      \tikzset{every left delimiter/.style={xshift=-1ex},every right delimiter/.style={xshift=1ex}}
      \matrix (mat) [matrix of math nodes, left delimiter=(, right delimiter=)] {
\draw[fill=gray] (0,0) rectangle (16,16);
\\
      };
\end{tikzpicture}
}
$\approx$
\resizebox{3.5cm}{!}{
\begin{tikzpicture}[baseline=-0.5ex]
      \tikzset{every left delimiter/.style={xshift=-1ex},every right delimiter/.style={xshift=1ex}}
      \matrix (mat) [matrix of math nodes, left delimiter=(, right delimiter=)] {
      \draw;
\foreach \i in {0,1,...,3} {
    \foreach \j in {0,4,...,12} {
        \draw[fill=gray] (16-\i-\j,\j) rectangle (16-\i-\j-1,\j+4);
    }
}
\\
      };
\end{tikzpicture}
}
\resizebox{3.5cm}{!}{
\begin{tikzpicture}[baseline=-0.5ex]
      \tikzset{every left delimiter/.style={xshift=-1ex},every right delimiter/.style={xshift=1ex}}
      \matrix (mat) [matrix of math nodes, left delimiter=(, right delimiter=)] {
      \draw;
\foreach \i in {0,1,...,3}{
    \foreach \j in {0,1,...,3}{
        \draw[fill=gray] (4*\i+\j,16-4*\j-\i) rectangle (4*\i+\j+1,15-4*\j-\i);
    }
}
\\
      };
\end{tikzpicture}
}
\resizebox{3.5cm}{!}{
\begin{tikzpicture}[baseline=-0.5ex]
      \tikzset{every left delimiter/.style={xshift=-1ex},every right delimiter/.style={xshift=1ex}}
      \matrix (mat) [matrix of math nodes, left delimiter=(, right delimiter=)] {
            \draw;
\foreach \i in {0,1,...,3} {
    \foreach \j in {0,4,...,12} {
        \draw[fill=gray] (\j,16-\i-\j) rectangle (\j+4,16-\i-\j-1);
    }
}
\\
      };
\end{tikzpicture}
}\\
\end{minipage}
\caption{The middle level factorization of
    a complementary low-rank matrix $K\approx U^2M^2(V^2)^*$
    where $N=n^2=4^2$ and $r=1$.
    Grey blocks indicate nonzero blocks.
    $U^2$ and $V^2$ are block-diagonal matrices with $4$ blocks.
    The diagonal blocks of $U^2$ and $V^2$ are assembled
    according to Equation \eqref{eq:expression-U} and \eqref{eq:expression-V}
    as indicated by gray rectangles.
    $M^2$ is a $4\times 4$ block matrix with each block
    $M^2_{\i,\j}$ itself being an $4\times 4$ block matrix
    containing diagonal weight matrix on the $(\j,\i)$ block.}
\label{fig:compression-2D}
\end{figure}
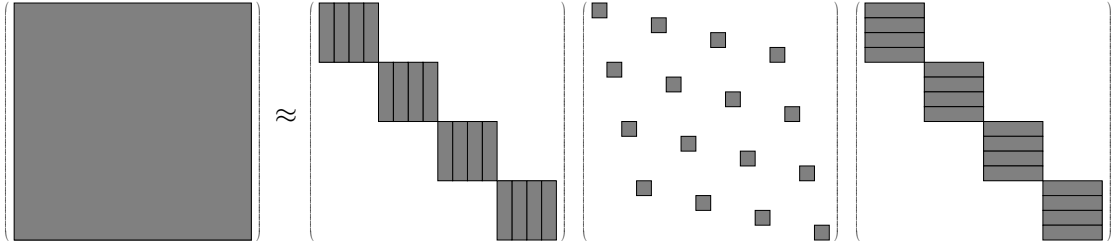

After computing the rank-$r$ factorization in \eqref{eq:Kij} for all
$\i$ and $\j$ in $\I^h$, we assemble all left factors $U^h_{\i,\j}$
into a matrix $U^h$, all middle factors into a matrix $M^h$, and all
right factors into a matrix $V^h$ so that
\begin{equation}
  \label{eq:UMV}
  K\approx U^hM^h(V^h)^*.
\end{equation}
Here $U^h$ is a block diagonal matrix of size $N\times rN$ with $n$
diagonal blocks $U^h_{\i}$ of size $n\times rn$:
\begin{equation*}
U^h=
\begin{pmatrix}
U^h_{{0}} & & &\\
 & U^h_{{1}} & &\\
 & & \ddots &\\
 & & & U^h_{{4^h-1}}
\end{pmatrix},
\end{equation*}
where each diagonal block $U^h_{\i}$ consists of the left factors
$U^h_{\i,\j}$ for all $\j$ as follows:
\begin{equation}
  U_{\i}^h=
  \begin{pmatrix}
    U^h_{\i,{0}} & U^h_{\i,{1}} & \cdots & U^h_{\i,{4^h-1}}
  \end{pmatrix}
  \in \bbC^{n \times rn}.
  \label{eq:expression-U}
\end{equation}
Similarly, $V^h$ is a block diagonal matrix of size $N\times rN$ with
$n$ diagonal blocks $V^h_{\j}$ of size $n\times rn$, where each
diagonal block $V^h_{\j}$ consists of the right factors $V^h_{\j,\i}$
for all $\i$ as follows:
\begin{equation}
  V^h_{\j}=
  \begin{pmatrix}
    V^h_{\j,{0}} & V^h_{\j,{1}} & \cdots & V^h_{\j,{4^h-1}}
  \end{pmatrix}
  \in \bbC^{n \times rn}.
  \label{eq:expression-V}
\end{equation}
The middle matrix $M^h\in \bbC^{rN\times rN}$ is an $n\times n$ block
matrix. The $(\i,\j)$-th block $M^h_{\i,\j}\in \bbC^{rn\times rn}$ is
itself an $n\times n$ block matrix. The only nonzero block of
$M^h_{\i,\j}$ is the $(\j,\i)$-th block, which is equal to the
$r\times r$ matrix $S^h_{\i,\j}$, and the other blocks of
$M^h_{\i,\j}$ are zero. We refer to Figure \ref{fig:compression-2D}
for a simple example of the middle level factorization when $N=4^2$.

\subsection{Recursive factorization}
\label{sec:rf}
In this section, we shall discuss how to recursively factorize
\begin{equation}
  \label{eq:RSU}
  U^\ell\approx U^{\ell+1}G^\ell
\end{equation}
and
\begin{equation}
  \label{eq:RSV}
  (V^\ell)^*\approx (H^\ell)^*(V^{\ell+1})^*
\end{equation}
for $\ell=h,h+1,\dots,L-1$. After these recursive factorizations, we
can construct the two-dimensional butterfly factorization
\begin{equation}
K\approx U^LG^{L-1}\cdots G^hM^h\left(H^h\right)^*\cdots
\left(H^{L-1}\right)^*\left(V^L\right)^*
\end{equation}
by substituting these recursive factorizations into \eqref{eq:UMV}.

\subsubsection{Recursive factorization of $U^h$}
\label{sec:rfu}

In the middle level factorization, we utilized the low-rank property
of $K^h_{\i,\j}$, the kernel matrix restricted in the domain
$A^h_\i\times B^h_\j\in T_X\times T_\Omega$, to obtain $U^h_{\i,\j}$
for $\i,\j\in\I^h$. We shall now use the complementary low-rank
property at level $\ell=h+1$, i.e., the matrix $K^{h+1}_{\i,\j}$
restricted in $A^{h+1}_\i\times B^{h-1}_\j\in T_X\times T_\Omega$ is
numerical low-rank for $\i\in\I^{h+1}$ and $\j\in\I^{h-1}$. These
factorizations of the column bases from level $h$ generate the column
bases at level $h+1$ through the following four steps: splitting,
merging, truncating, and assembling.

\paragraph{Splitting.} In the middle level factorization, we have constructed
\begin{equation*}
  U^h=
  \begin{pmatrix}
    U^h_{{0}} & & &\\
    & U^h_{{1}} & &\\
    & & \ddots &\\
    & & & U^h_{{4^h-1}}
  \end{pmatrix}
  \quad
  \text{with}
  \quad
  U_{\i}^h=
  \begin{pmatrix}
    U^h_{\i,{0}} & U^h_{\i,{1}} & \cdots & U^h_{\i,{4^h-1}}
  \end{pmatrix}\in \bbC^{n\times r n},
\end{equation*}
where each $U^h_{\i,\j}\in \bbC^{n\times r}$. Each node $A^h_\i$ in the
quadtree $T_X$ on the level $h$ has four child nodes on the level
$h+1$, denoted by $\{A^{h+1}_{4\i+{t}}\}_{t=0,1,2,3}$.  According to
this structure, one can split $U^h_{\i,\j}$ into four parts in the row
space,
\begin{equation}\label{eq:splitandmerge}
  U^h_{\i,\j}=
  \begin{pmatrix}
    U^{h,0}_{\i,\j}\\
    \midrule
    U^{h,1}_{\i,\j}\\
    \midrule
    U^{h,2}_{\i,\j}\\
    \midrule
    U^{h,3}_{\i,\j}
  \end{pmatrix},
\end{equation}
where $U^{h,t}_{\i,\j}$ approximately spans the column space of the
submatrix of $K$ restricted to $A^{h+1}_{4\i+{t}}\times B^h_{\j}$ for
each $t=0,\ldots,3$. Combining this with the definition of $U^h_i$
gives rise to
\begin{equation}
  U^h_{\i}=
  \begin{pmatrix}
    U^h_{\i,{0}} & U^h_{\i,{1}} & \cdots & U^h_{\i,{4^h-1}}
  \end{pmatrix}
  =
  \begin{pmatrix}
    U^{h,0}_{\i,{0}} & U^{h,0}_{\i,{1}} &\dots
    & U^{h,0}_{\i,{4^h-1}}\\
    \midrule
    U^{h,1}_{\i,{0}} & U^{h,1}_{\i,{1}} &\dots
    & U^{h,1}_{\i,{4^h-1}}\\
    \midrule
    U^{h,2}_{\i,{0}} & U^{h,2}_{\i,{1}} &\dots
    & U^{h,2}_{\i,{4^h-1}}\\
    \midrule
    U^{h,3}_{\i,{0}} & U^{h,3}_{\i,{1}} &\dots
    & U^{h,3}_{\i,{4^h-1}}
  \end{pmatrix}
  =:
  \begin{pmatrix}
    U^{h,0}_{\i}\\
    \midrule
    U^{h,1}_{\i}\\
    \midrule
    U^{h,2}_{\i}\\
    \midrule
    U^{h,3}_{\i}
  \end{pmatrix},
\end{equation}
where $U^{h,t}_{\i}$ approximately spans the column space of the
matrix $K$ restricted to $A^{h+1}_{4\i+{t}}\times \Omega$.

\paragraph{Merging.} 
The merging step merges adjacent matrices $U^{h,t}_{\i,\j}$ in the
column space to obtain low-rank matrices. For any $\i\in \I^h$ and
$\j\in \I^{h-1}$, the merged matrix
\begin{equation}
  \begin{pmatrix}
    U^{h,t}_{\i,4\j+{0}} &  U^{h,t}_{\i,4\j+{1}} & U^{h,t}_{\i,4\j+{2}} &  U^{h,t}_{\i,4\j+{3}}
  \end{pmatrix}
  \in \bbC^{n/4\times 4r}
  \label{eq:blocks}
\end{equation}
approximately spans the column space of $K^{h+1}_{4\i+{t},\j}$
corresponding to the domain $A_{4\i+t}^{h+1}\times B^{h-1}_{\j}$. By
the complementary low-rank property of the matrix $K$, we know
$K^{h+1}_{4\i+{t},\j}$ is numerically low-rank. Hence, the matrix in
\eqref{eq:blocks} is also a numerically low-rank matrix.  This is the
merging step equivalent to moving from level $h$ to level $h-1$ in
$T_\Omega$.

\paragraph{Truncating.}
The third step computes its rank-$r$ approximation using the standard
truncated {SVD} and putting it to the form of \eqref{eq:lowrank2}. For
each $\i\in\I^h$ and $\j\in\I^{h-1}$, the factorization
\begin{equation}
  \begin{pmatrix}
    U^{h,t}_{\i,4\j+{0}} &  U^{h,t}_{\i,4\j+{1}} &    U^{h,t}_{\i,4\j+{2}} &  U^{h,t}_{\i,4\j+{3}}
  \end{pmatrix}
  \approx U^{h+1}_{4\i+{t},\j}G^{h}_{4\i+{t},\j},
  \label{eq:blockF}
\end{equation}
defines $U^{h+1}_{4\i+{t},\j}\in \bbC^{n/4\times r}$ and
$G^{h}_{4\i+{t},\j}\in\bbC^{r\times 4r}$.

\paragraph{Assembling} In the final step, we construct the factorization 
$U^h \approx U^{h+1} G^h$ using \eqref{eq:blockF}. Since 
$\I^{h+1}$ is the same as $\{4\i+\ve{t}\}_{\i\in\I^h,t=0,1,2,3}$,
one can arrange \eqref{eq:blockF} for all $\i$ and $\j$ into a single
formula as follows:
\begin{equation*}
  \begin{split}
    &U^h \approx U^{h+1}G^{h} =\\
    &\begin{pmatrix}
       U^{h+1}_{{0}} & & & & & & & & &\\
       & \ddots & & & & & & & &\\
       & & U^{h+1}_{{3}} & & & & & & &\\
       & & & U^{h+1}_{{4}} & & & & &\\
       & & & & \ddots & & & &\\
       & & & & & U^{h+1}_{{7}} & & &\\
       & & & & & & \ddots & &\\
       & & & & & & & U^{h+1}_{{4^{h+1}-4}} &\\
       & & & & & & & & \ddots \\
       & & & & & & & & & U^{h+1}_{{4^{h+1}-1}}
     \end{pmatrix}
    \begin{pmatrix}
      G^{h}_{{0}} & & &\\
      \vdots & & &\\
      G^{h}_{{3}} & & &\\
      & G^{h}_{{4}} & &\\
      & \vdots & &\\
      & G^{h}_{{7}} & &\\
      & & \ddots &\\
      & & & G^{h}_{{4^{h+1}-4}}\\
      & & & \vdots\\
      & & & G^{h}_{{4^{h+1}-1}}
    \end{pmatrix},
  \end{split}
\end{equation*}
where the blocks are given by
\begin{equation*}
  U_{\i}^{h+1}=
  \begin{pmatrix}
    U^{h+1}_{\i,{0}} & U^{h+1}_{\i,{1}} & \cdots &
    U^{h+1}_{\i,{4^{h-1}-1}}
  \end{pmatrix}
\end{equation*}
and
\begin{equation*}
  G^{h}_{\i}=
  \begin{pmatrix}
    G^{h}_{\i,{0}} & & &\\
    & G^{h}_{\i,{1}} & &\\
    & & \ddots &\\
    & & & G^{h}_{\i,{4^{h-1}-1}}\\
  \end{pmatrix}
\end{equation*}
for $\i\in\I^{h+1}$. Figure \ref{fig:compression-2DU} shows a toy
example of the recursive factorization of $U^h$ when $N = 4^2$, $h=2$
and $r=1$. Since there are $\O(1)$ nonzero entries in each
$G^h_{\i,\j}$ and $\O(4^{h+1}\cdot 4^{h-1})=\O(N)$ such matrices,
there are only $\O(N)$ nonzero entries in $G^h$.

\begin{figure}[htp]
\begin{minipage}{\textwidth}
\centering
\resizebox{3.5cm}{!}{
\begin{tikzpicture}[baseline=-0.5ex]
      \tikzset{every left delimiter/.style={xshift=-1ex},every right delimiter/.style={xshift=1ex}}
      \matrix (mat) [matrix of math nodes, left delimiter=(, right delimiter=)] {
\draw;
\foreach \i in {0,1,...,3} {
    \foreach \j in {0,4,...,12} {
        \draw[fill=gray] (16-\i-\j,\j) rectangle (16-\i-\j-1,\j+1);
        \draw[fill=gray] (16-\i-\j,\j+1) rectangle (16-\i-\j-1,\j+2);
        \draw[fill=gray] (16-\i-\j,\j+2) rectangle (16-\i-\j-1,\j+3);
        \draw[fill=gray] (16-\i-\j,\j+3) rectangle (16-\i-\j-1,\j+4);
    }
}
\\
      };
\end{tikzpicture}
}
$\approx$
\resizebox{3.5cm}{!}{
\begin{tikzpicture}[baseline=-0.5ex]
      \tikzset{every left delimiter/.style={xshift=-1ex},every right delimiter/.style={xshift=1ex}}
      \matrix (mat) [matrix of math nodes, left delimiter=(, right delimiter=)] {
\draw;
\foreach \k in {0,4,...,12} {
\foreach \i in {0,...,3} {
        \draw[fill=gray] (\k,12-\k+\i) rectangle (\k+1,12-\k+\i+1);
        \draw[fill=lightgray] (\k+1,12-\k+\i) rectangle (\k+4,12-\k+\i+1);
}
}
\\
      };
\end{tikzpicture}
}
=
\resizebox{3.5cm}{!}{
\begin{tikzpicture}[baseline=-0.5ex]
      \tikzset{every left delimiter/.style={xshift=-1ex},every right delimiter/.style={xshift=1ex}}
      \matrix (mat) [matrix of math nodes, left delimiter=(, right delimiter=)] {
\draw;
\foreach \k in {0,1,...,15}{
    \draw[fill=gray] (16-\k,\k) rectangle (16-\k-1,\k+1);
}\\
      };
\end{tikzpicture}
}
\resizebox{3.5cm}{!}{
\begin{tikzpicture}[baseline=-0.5ex]
      \tikzset{every left delimiter/.style={xshift=-1ex},every right delimiter/.style={xshift=1ex}}
      \matrix (mat) [matrix of math nodes, left delimiter=(, right delimiter=)] {
\draw;
\foreach \k in {0,4,...,12}{
\foreach \i in {0,1,...,3} {
    \draw[fill=lightgray] (16-\k,\k+\i) rectangle (12-\k,\k+\i+1);
}
}\\
      };
\end{tikzpicture}
}\\
\resizebox{3.5cm}{!}{
}
\end{minipage}
\caption{The recursive factorization of $U^2$
    in Figure~\ref{fig:compression-2D}.
    Left matrix: $U^2$ with each diagonal block
    partitioned into smaller blocks according to
    Equation \eqref{eq:splitandmerge} as indicated by black rectangles;
    Middle-left matrix: low-rank approximations of submatrices
    in $U^2$ given by Equation \eqref{eq:blockF};
    Middle right matrix: $U^3$;
    Right matrix: $G^2$.}
\label{fig:compression-2DU}
\end{figure}

In a similar way, we can now factorize $U^\ell\approx
U^{\ell+1}G^\ell$ for $h< \ell\leq L-1$. As before, the key point is
that the columns of 
\begin{equation}
  \begin{pmatrix}
    U^{\ell,t}_{\i,4\j+{0}} & U^{\ell,t}_{\i,4\j+{1}} &
    U^{\ell,t}_{\i,4\j+{2}} & U^{\ell,t}_{\i,4\j+{3}}
  \end{pmatrix}
\end{equation}
approximately span the column space of $K^{\ell+1}_{4\i+{t},\j}$,
which is of rank $r$ numerically due to the complementary low-rank
property.  Computing its rank-$r$ approximation via the standard
truncated {SVD} results in a form of \eqref{eq:lowrank2}
\begin{equation}
  \begin{pmatrix}
    U^{\ell,t}_{\i,4\j+{0}} & U^{\ell,t}_{\i,4\j+{1}} &
    U^{\ell,t}_{\i,4\j+{2}} & U^{\ell,t}_{\i,4\j+{3}}
  \end{pmatrix}
  \approx U^{\ell+1}_{4\i+{t},\j} G^{\ell}_{4\i+{t},\j}
  \label{eq:RSh2}
\end{equation}
for $\i\in\I^\ell$ and $\j\in \I^{L-\ell-1}$. After
assembling these factorizations together, we obtain
\begin{equation*}
  \begin{split}
    &U^\ell \approx U^{\ell+1}G^{\ell} =\\
    &\begin{pmatrix}
       U^{\ell+1}_{{0}} & & & & & & & & &\\
       & \ddots & & & & & & & &\\
       & & U^{\ell+1}_{{3}} & & & & & & &\\
       & & & U^{\ell+1}_{{4}} & & & & &\\
       & & & & \ddots & & & &\\
       & & & & & U^{\ell+1}_{{7}} & & &\\
       & & & & & & \ddots & &\\
       & & & & & & & U^{\ell+1}_{{4^{\ell+1}-4}} &\\
       & & & & & & & & \ddots \\
       & & & & & & & & & U^{\ell+1}_{{4^{\ell+1}-1}}
     \end{pmatrix}
    \begin{pmatrix}
      G^{\ell}_{{0}} & & &\\
      \vdots & & &\\
      G^{\ell}_{{3}} & & &\\
      & G^{\ell}_{{4}} & &\\
      & \vdots & &\\
      & G^{\ell}_{{7}} & &\\
      & & \ddots &\\
      & & & G^{\ell}_{{4^{\ell+1}-4}}\\
      & & & \vdots\\
      & & & G^{\ell}_{{4^{\ell+1}-1}}
    \end{pmatrix},
  \end{split}
\end{equation*}
where
\begin{equation*}
  U_{\i}^{\ell+1}=
  \begin{pmatrix}
    U^{\ell+1}_{\i,{0}} & U^{\ell+1}_{\i,{1}} & \cdots &
    U^{\ell+1}_{\i,{4^{L-\ell-1}-1}}
  \end{pmatrix}
\end{equation*}
and
\begin{equation*}
  G^{\ell}_{\i}=
  \begin{pmatrix}
    G^{\ell}_{\i,{0}} & & &\\
    & G^{\ell}_{\i,{1}} & &\\
    & & \ddots &\\
    & & & G^{\ell}_{\i,{4^{L-\ell-1}-1}}\\
\end{pmatrix}
\end{equation*}
for $\i \in \I^{\ell+1}$.

After the $L-h$ step of recursive factorizations $U^\ell \approx
U^{\ell+1}G^{\ell}$ for $\ell=h,h+1,\dots,L-1$, the recursive
factorization of $U^h$ takes the following form:
\begin{equation}
  \label{eq:RFUh}
  U^h\approx U^{L}G^{L-1}\cdots G^{h}.
\end{equation}
Similarly to the analysis of $G^{h}$, it is also easy to check that
there are only $\O(N)$ nonzero entries in each $G^{\ell}$ in
\eqref{eq:RFUh}. As to the first factor $U^{L}$, it has
$\O(N)$ nonzero entries since there are $\O(N)$ diagonal blocks in
$U^L$ and each block contains $\O(1)$ entries.

\subsubsection{Recursive factorization of $V^h$}
\label{sec:rfv}

The recursive factorization of $V^\ell$ is similar to that of $U^\ell$
for $\ell=h,h+1,\dots,L-1$.  At each level $\ell$, we benefit from the
fact that
\[
\begin{pmatrix}
  V^{\ell,t}_{\j,4\i+{0}} &
  V^{\ell,t}_{\j,4\i+{1}} &
  V^{\ell,t}_{\j,4\i+{2}} &
  V^{\ell,t}_{\j,4\i+{3}}
\end{pmatrix}
\]
approximately spans the row space of $K^{L-\ell-1}_{\i,4\j+{t}}$ and
hence is numerically low-rank for $\j\in\I^{L-\ell}$ and
$\i\in\I^{\ell-1}$.  Applying the same procedure in Section
\ref{sec:rfu} to $V^{h}$ leads to
\begin{equation}
  \label{eq:VF}
  V^h\approx V^LH^{L-1}\cdots H^h.
\end{equation}


\subsection{Complexity analysis}
\label{sec:ca}

\begin{figure}[htp]
\begin{minipage}{\textwidth}
\centering
\vcenteredinclude{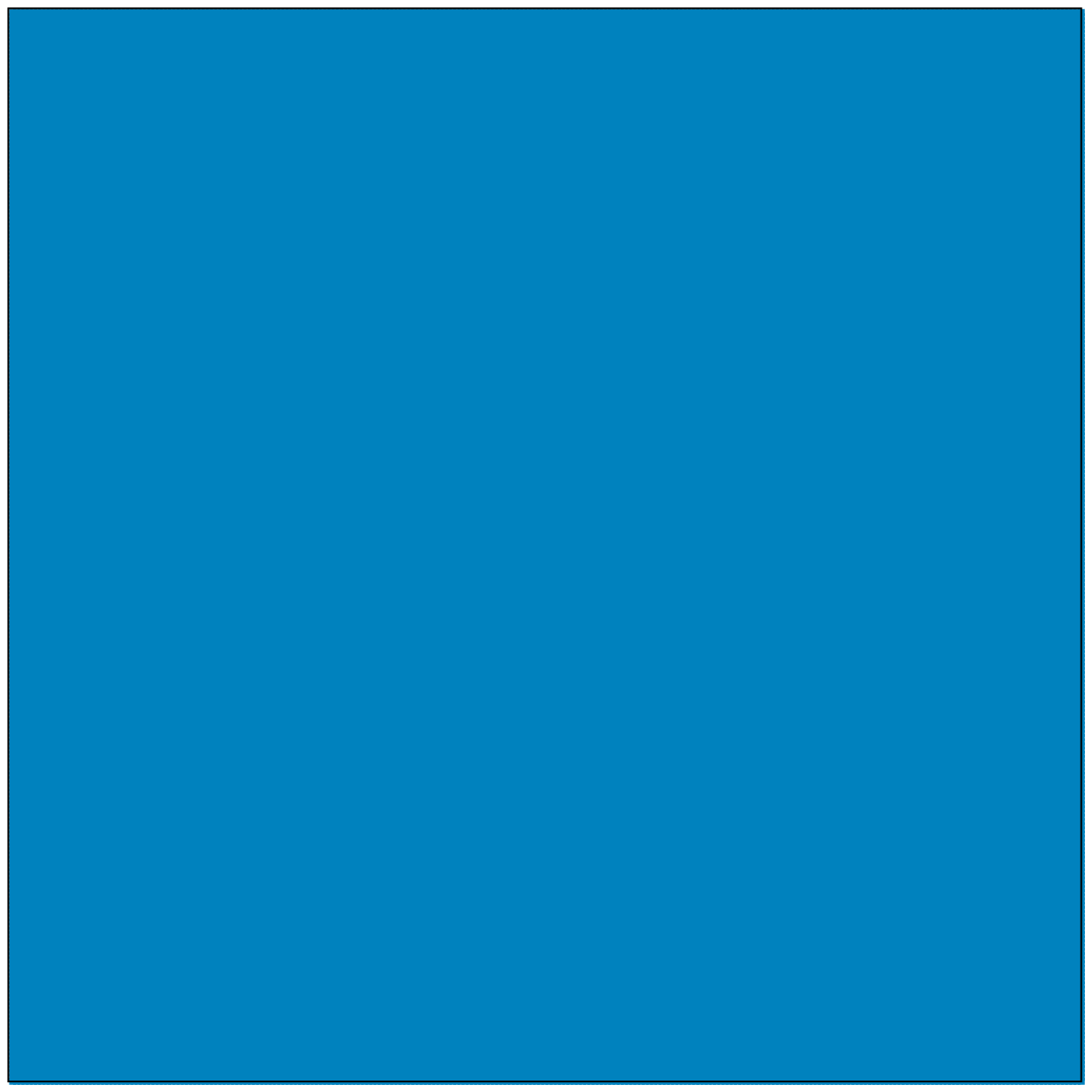}
$\approx$
\vcenteredinclude{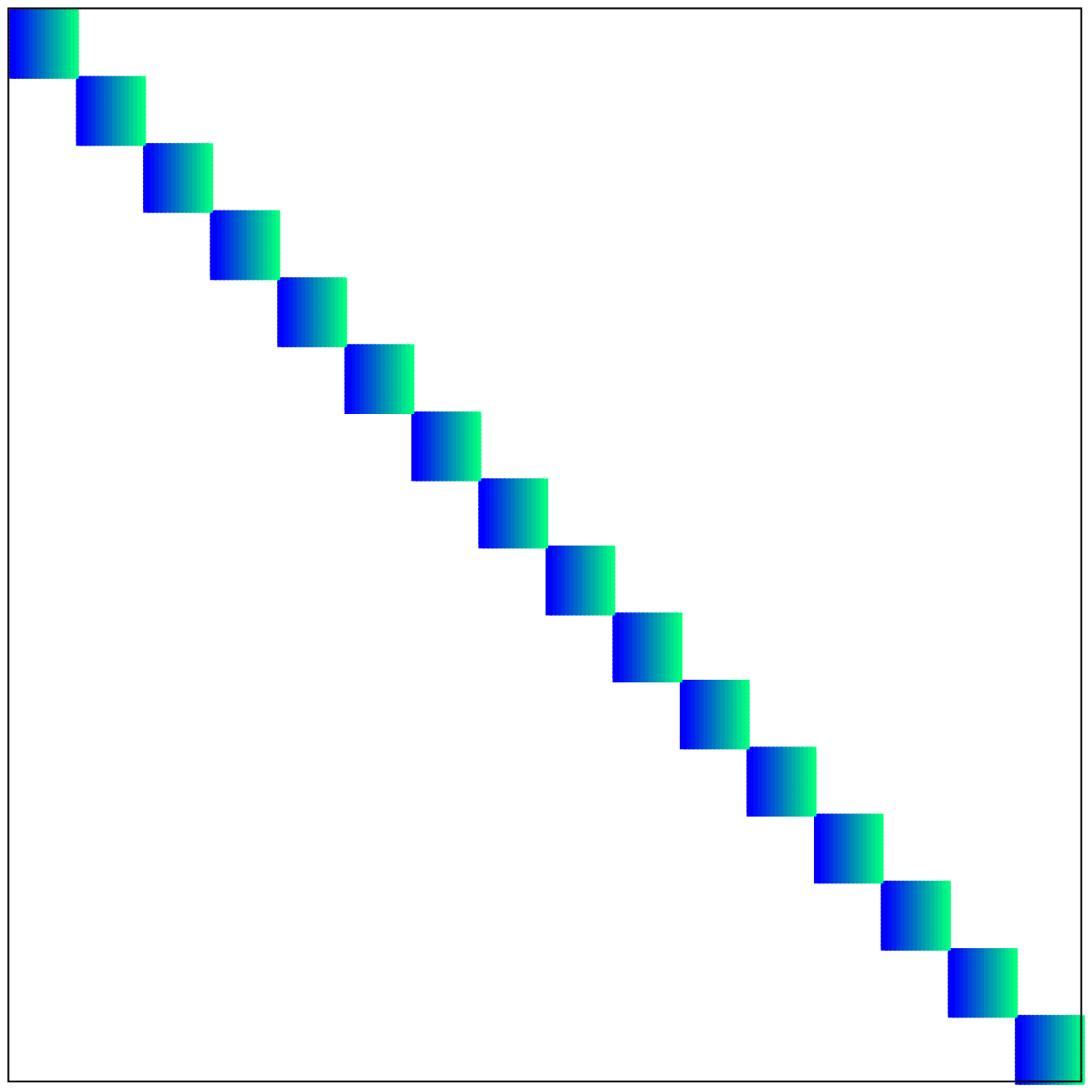}
\vcenteredinclude{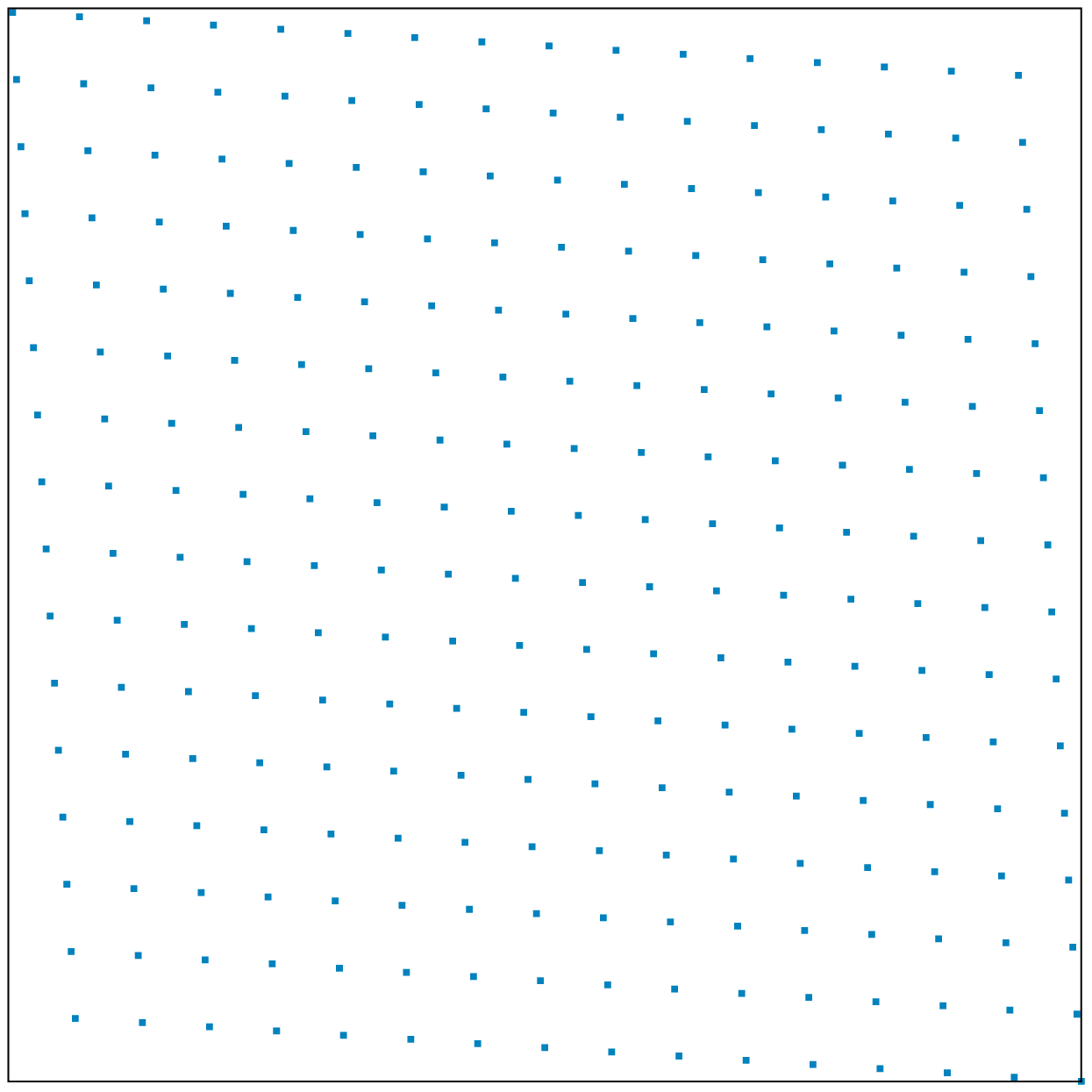}
\vcenteredinclude{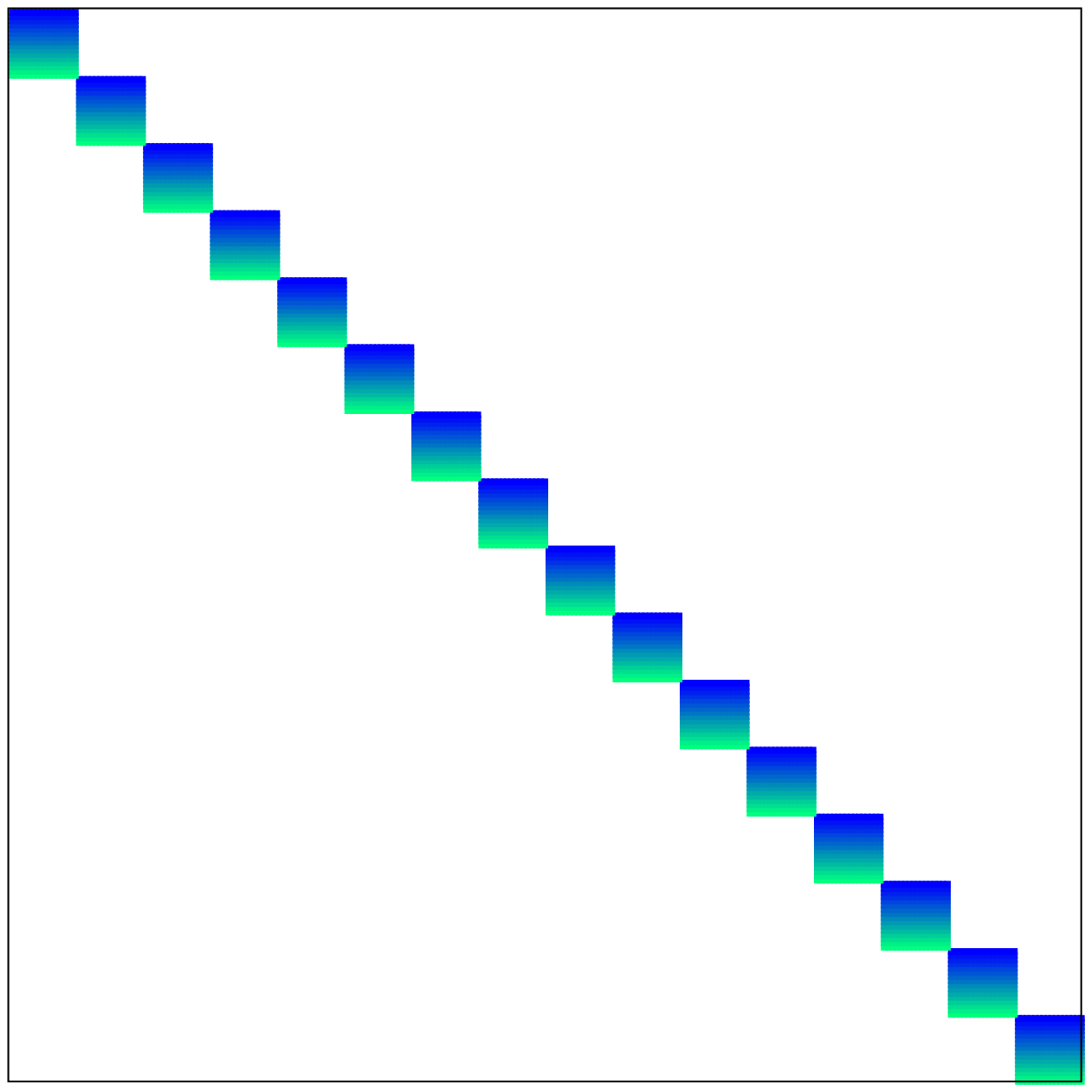}\\
\vspace{0.35cm}
$\approx$
\vcenteredinclude{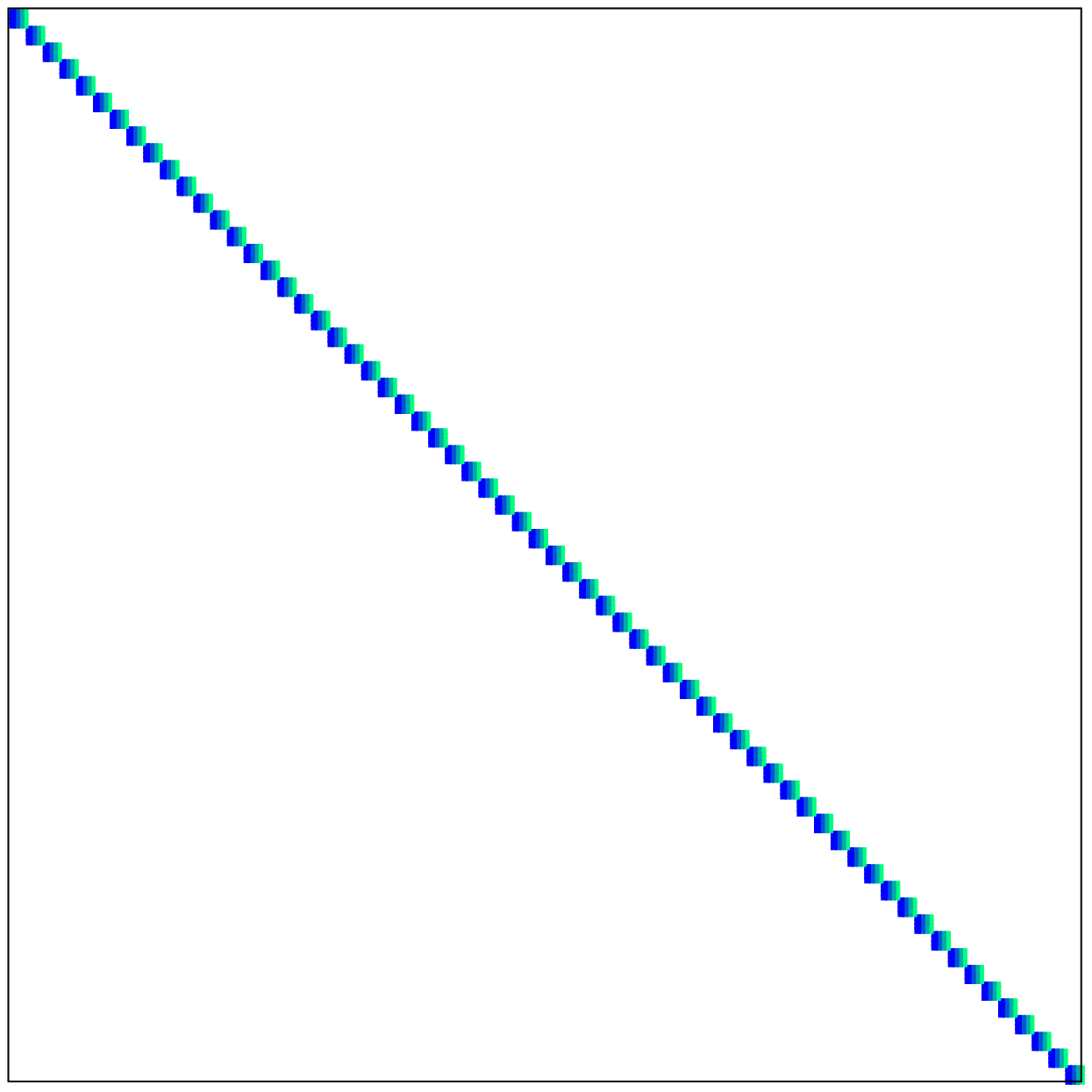}
\vcenteredinclude{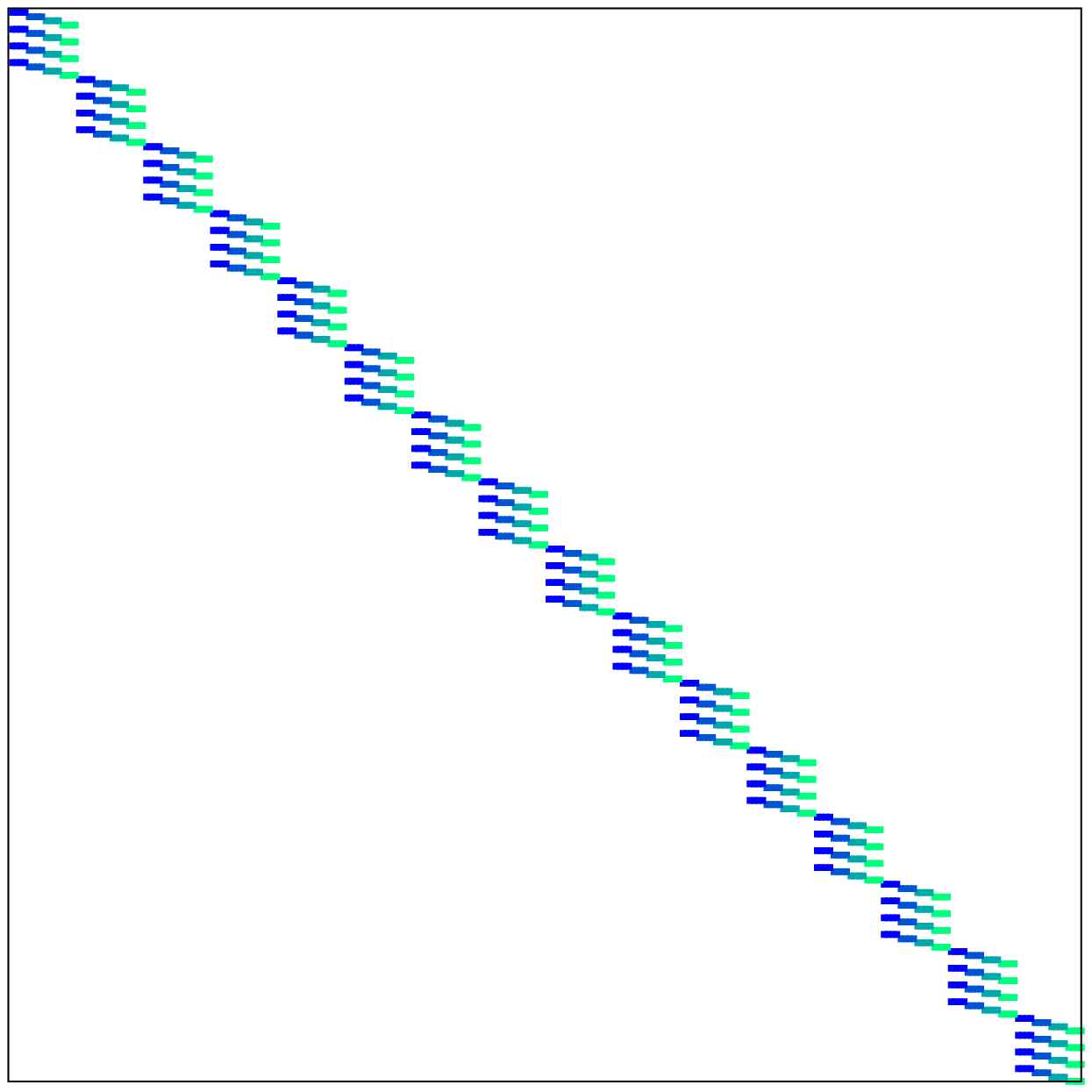}
\vcenteredinclude{figure/fig-comp-modules-real/fig-comp-M}\\
\vspace{0.35cm}
\vcenteredinclude{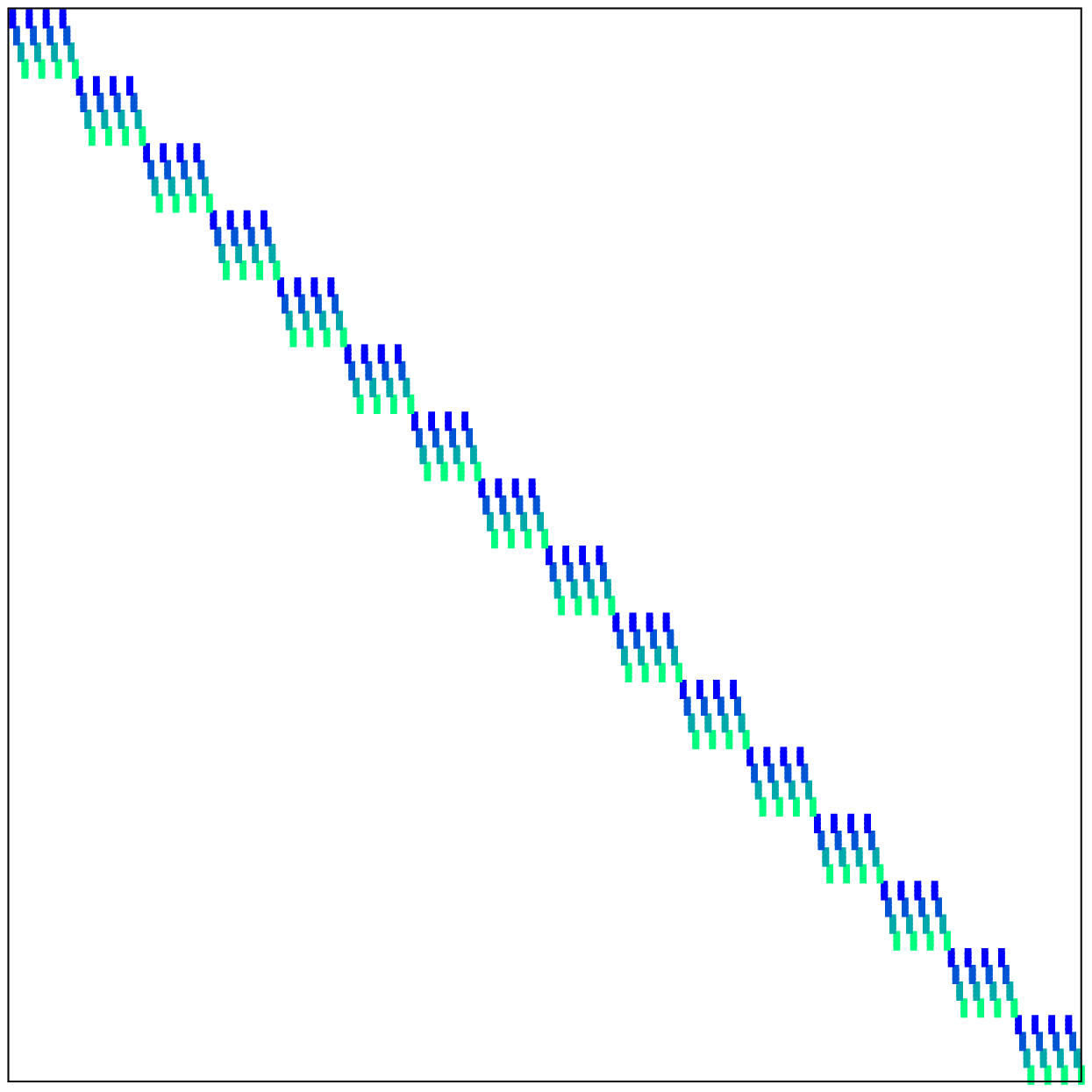}
\vcenteredinclude{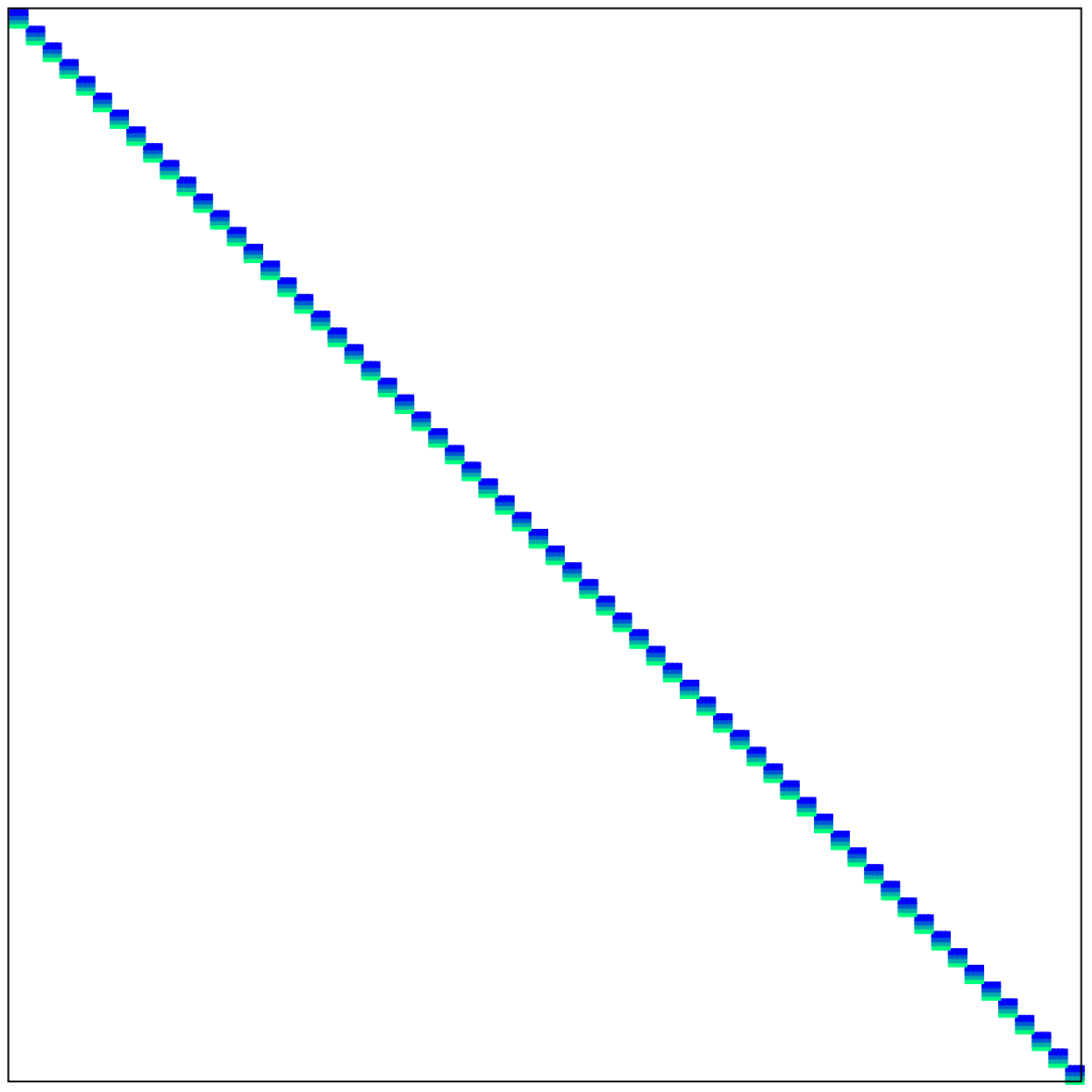}\\
\vspace{0.35cm}
$\approx$
\vcenteredinclude{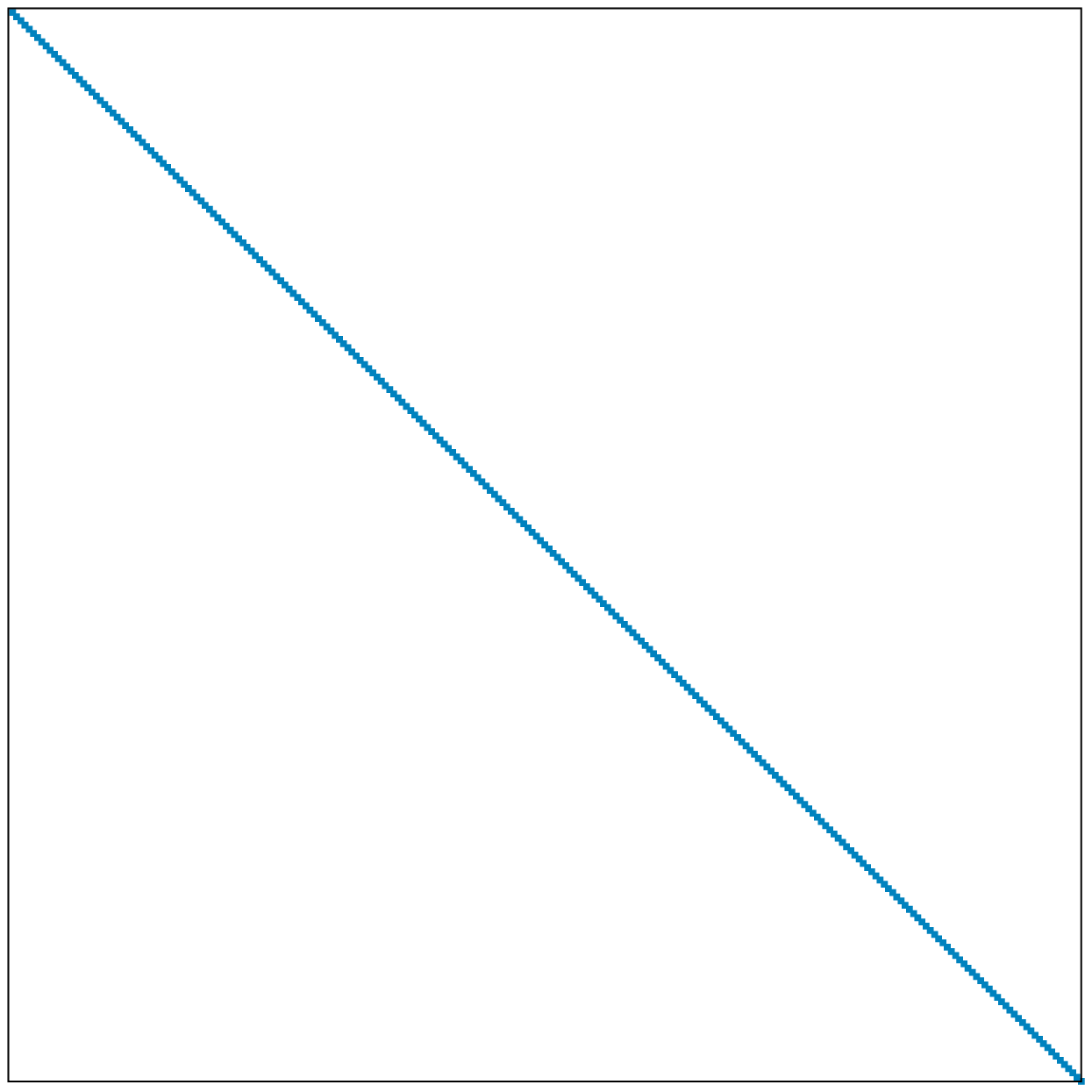}
\vcenteredinclude{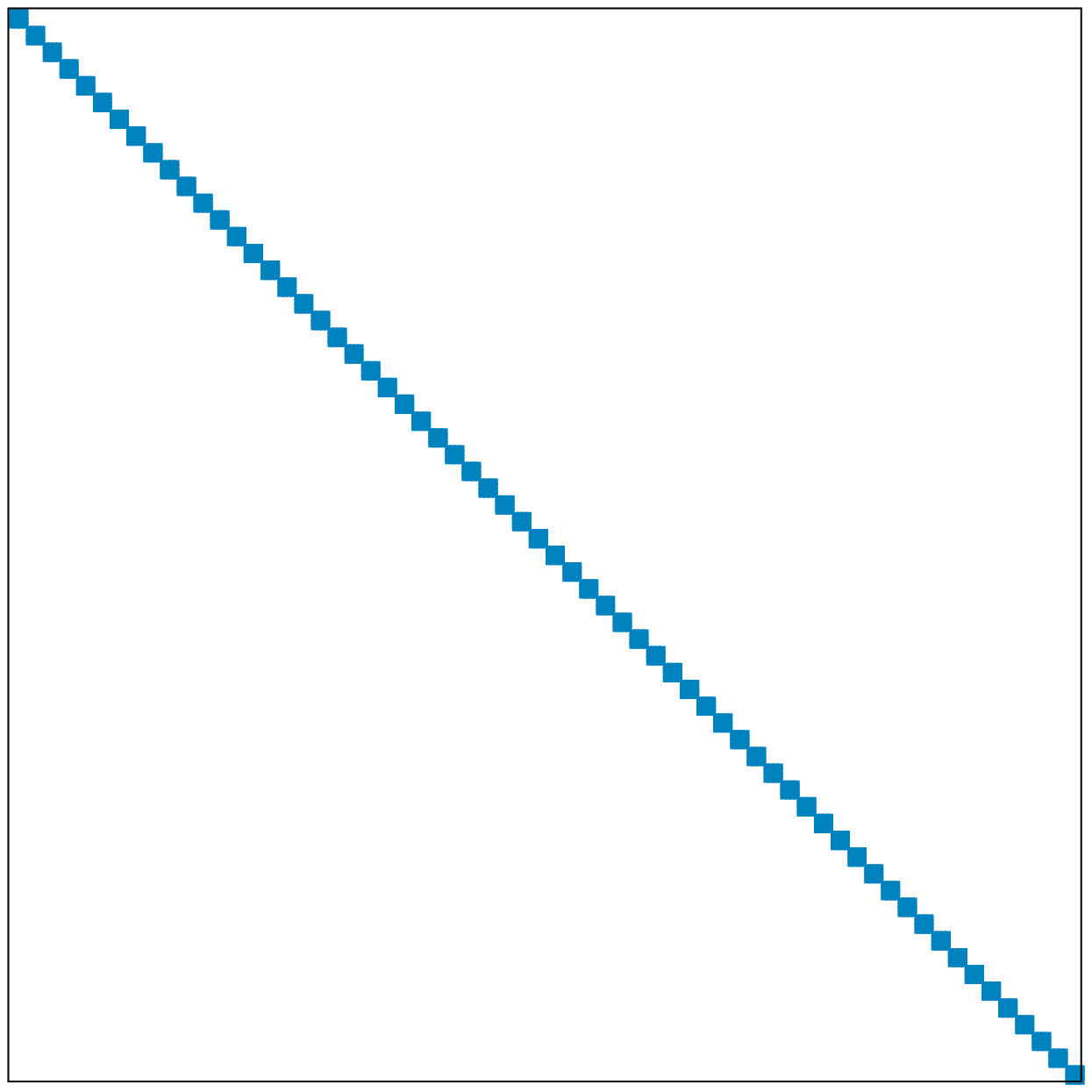}
\vcenteredinclude{figure/fig-comp-modules-real/fig-comp-G0}
\vcenteredinclude{figure/fig-comp-modules-real/fig-comp-M}\\
\vspace{0.35cm}
\vcenteredinclude{figure/fig-comp-modules-real/fig-comp-H0}
\vcenteredinclude{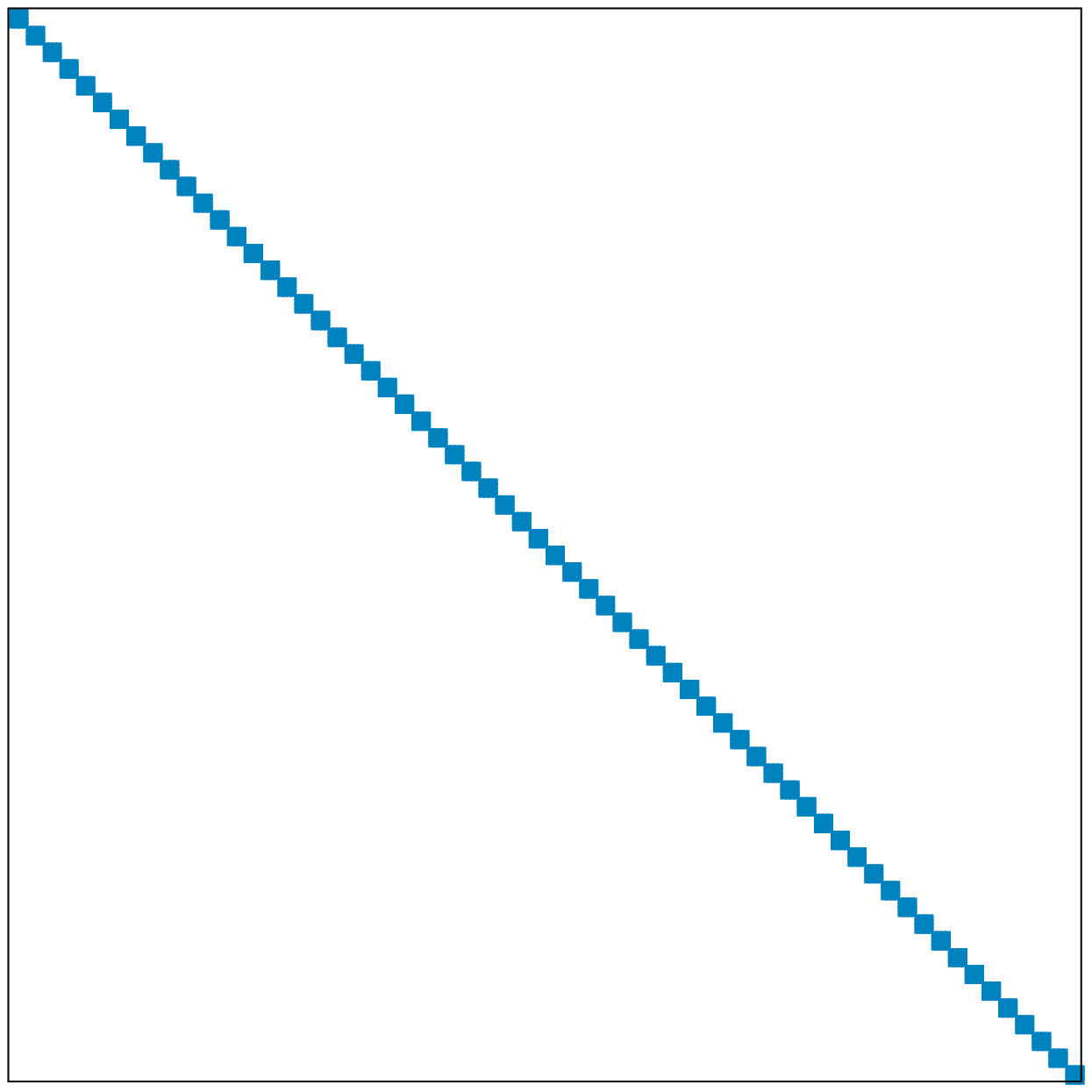}
\vcenteredinclude{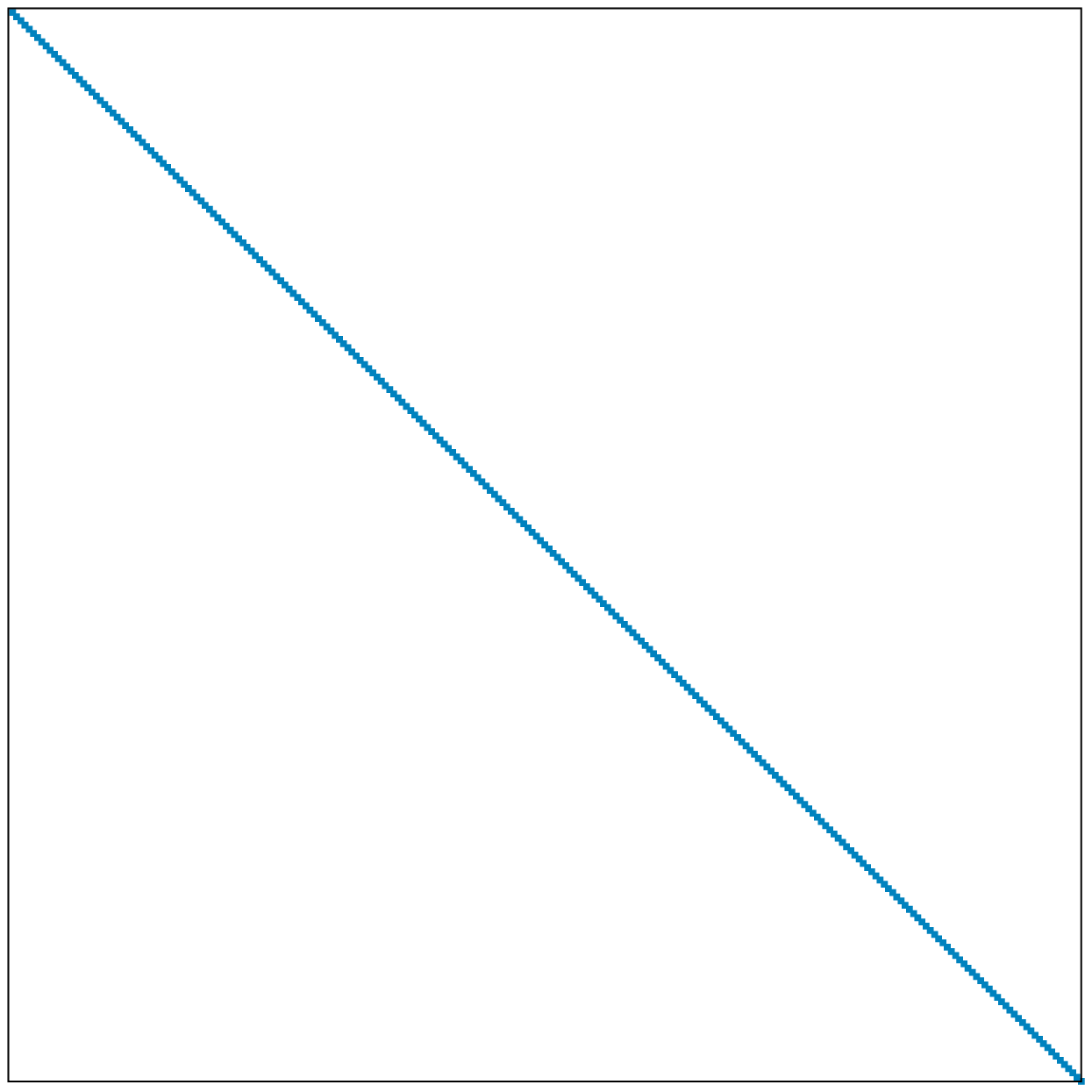}\\
\end{minipage}
\caption{A full butterfly factorization for a two dimensional problem
    of size $16^2$ and fixed rank $r=1$.
    The above figure visualizes the matrices in $K \approx U^3M^3(V^3)^*
    \approx U^4G^3M^3(H^3)^*(V^4)^* \approx U^5G^4G^3M^3(H^3)^*(H^4)^*(V^5)^*$.
}
\label{fig:compression-2D-real}
\end{figure}

By combining the results of the middle level factorization in
\eqref{eq:UMV} and the recursive factorizations in \eqref{eq:RFUh} and
\eqref{eq:VF}, we obtain the final butterfly factorization
\begin{equation}
K\approx U^LG^{L-1}\cdots G^hM^h\left(H^h\right)^*\cdots
\left(H^{L-1}\right)^*\left(V^L\right)^*,
\end{equation}
each factor of which contains $\O(N)$ nonzero entries.  We refer to
Figure \ref{fig:compression-2D-real} for an illustration of the
butterfly factorization of $K$ when $N=16^2$.

The complexity of constructing the butterfly factorization comes from two
parts: the middle level factorization and the recursive factorization.
For the middle level factorization, the construction cost is different
depending on which of the two cases mentioned in Section
\ref{sec:mlf} is under consideration, since they use different
approaches in constructing rank-$r$ SVDs at the middle level.
\begin{itemize}
\item In Case (i), the dominant cost is to apply $K$ and $K^*$ to
  $N^{1/2}$ Gaussian random matrices of size $N\times O(1)$.  Assuming
  that the given black-box routine for applying $K$ and $K^*$ to a
  vector takes $\O(C_K(N))$ operations, the total operation complexity
  is $\O(C_K(N)N^{1/2})$.
\item In Case (ii), we apply the SVD procedure with random sampling to
  $N$ submatrices of size $N^{1/2}\times N^{1/2}$.  Since the
  operation complexity for each submatrix is $\O(N^{1/2})$, the
  overall complexity is $\O(N^{3/2})$.
\end{itemize}

In the recursive factorization stage, most of the work comes from
factorizing $U^h$ and $V^h$. There are $\O(\log N)$ stages appeared in
the factorization of $U^h$. At the $\ell$ stage, the matrix $U^\ell$
to be factorized consists of $4^\ell$ diagonal blocks.  There are
$\O(N)$ factorizations and each factorization takes $\O(N/4^\ell)$
operations.  Hence, the operation complexity to factorize $U^\ell$ is
$\O(N^2/4^\ell)$.  Summing up all the operations in each step yields
the overall operation complexity for recursively factorizing $U^h$:
\begin{equation}
\sum^{L-1}_{\ell=h}\O(N^2/4^\ell) = \O(N^{3/2}).
\end{equation}

The peak of the memory usage of the butterfly factorization is due
to the middle level factorization where we need to store the results
of $\O(N)$ factorizations of size $\O(N^{1/2})$. Hence, the memory
complexity for the two-dimensional butterfly factorization is
$\O(N^{3/2})$. For Case (ii), one can actually do better by following
the same argument in \cite{1dbf}. One can interleave the order of
generation and recursive factorization of $U^h_{\i,\j}$ and
$V^h_{\j,\i}$. By factorizing $U^h_{\i,\j}$ and $V^h_{\j,\i}$
individually instead of formulating \eqref{eq:UMV}, the memory
complexity in Case (ii) can be reduced to $\O(N\log N)$.

The cost of applying the butterfly factorization is equal to the
number of nonzero entries in the final factorization, which is
$\O(N\log N)$. Table \ref{tab:complexity} summarizes the complexity
analysis for the two-dimensional butterfly factorization.

\begin{table}[ht!]
  \centering
  \begin{tabular}{llcc}
    \toprule
    & & SVD via rand. matvec & SVD via rand. sampling \\
    \toprule
    \multirow{5}{3cm}{Factorization Complexity} &    \parbox{2.5cm}{Middle level\\factorization}
    & $\O(C_K(N) N^{1/2} )$ & $\O(N^{3/2})$ \\
    \cmidrule(r){2-4}
    & \parbox{2.5cm}{Recursive\\factorization} & \multicolumn{2}{c}{$\O(N^{3/2})$}\\
    \cmidrule(r){2-4}
    & Total & $\O(C_K(N) N^{1/2} )$ & $\O(N^{3/2})$\\
    \toprule
    \parbox{3cm}{Memory\\Complexity} & & $\O(N^{3/2})$ & $\O(N\log N)$ \\
    \toprule
    \parbox{3cm}{Application\\Complexity} & & \multicolumn{2}{c}{$\O(N\log N)$}\\
    \bottomrule
  \end{tabular}
  \caption{The time and memory complexity of the two-dimensional
    butterfly factorization. Here $C_K(N)$ is the complexity of
    applying the matrices $K$ and $K^*$ to a vector. For most
    butterfly algorithms, $C_K(N)=O(N\log N)$.
  }
  \label{tab:complexity}
\end{table}

\subsection{Extensions}\label{sec:ext}

We have introduced the two-dimensional butterfly factorization for a
complementary low-rank kernel matrix $K$ in the entire domain $X\times
\Omega$.  Although we have assumed the uniform grid in \eqref{eq:X}
and \eqref{eq:Omega}, the butterfly factorization extends naturally to
more general settings.

In the case with non-uniform point sets $X$ or $\Omega$, one can still
construct a butterfly factorization for $K$ following the same
procedure.  More specifically, we still construct two trees $T_X$ and
$T_\Omega$ {\em adaptively} via hierarchically partitioning the square
domains covering $X$ and $\Omega$.  For non-uniform point sets $X$ and
$\Omega$, the numbers of points in $A^\ell_\i$ and $B^{L-\ell}_\j$ are
different. If a node does not contain any point inside it, it is
simply discarded from the quadtree.

The complexity analysis summarized in Table \ref{tab:complexity}
remains valid in the case of non-uniform point sets $X$ and $\Omega$.
On each level $\ell = h,\dots,L$ of the butterfly factorization,
although the sizes of low-rank submatrices are different, the total
number of submatrices and the numerical rank remain the same. Hence,
the total operation and memory complexity remains the same as
summarized in Table \ref{tab:complexity}.


\section{Polar Butterfly Factorization}
\label{sec:pbf}
In Section \ref{sec:gbf}, we have introduced a two-dimensional
{butterfly factorization} for a complementary low-rank kernel matrix
$K$ in the entire domain $X\times \Omega$.  In this section, we will
introduce a polar butterfly factorization to deal with the kernel
function $K(x,\xi)=e^{2\pi \imath \Phi(x,\xi)}$. Such a kernel matrix
has a singularity at $\xi=0$ and the approach taken here follows the
polar butterfly algorithm proposed in \cite{fio09}.

\subsection{Polar butterfly algorithm}
\label{sec:pba}

The multidimensional {Fourier} integral operator (FIO) is defined as
\begin{equation}\label{eq:fio}
  u(x) = \sum_{\xi\in\Omega} e^{2\pi \imath \Phi(x,\xi)}g(\xi),
  \quad x\in X,
\end{equation}
where the phase function $\Phi(x,\xi)$ is assumed to be real-analytic
in $(x,\xi)$ for $\xi\neq 0$, and is homogeneous of degree 1 in $\xi$,
namely, $\Phi(x,\lambda\xi) = \lambda\Phi(x,\xi)$ for all $\lambda>0$.
Here the grids $X$ and $\Omega$ are the same as those in
\eqref{eq:X} and \eqref{eq:Omega}.

As the phase function $\Phi(x,\xi)$ is singular at $\xi=0$, the
numerical rank of the kernel $e^{2\pi \imath \Phi(x,\xi)}$ in a domain
near or containing $\xi=0$ is typically large. Hence, in general
$K(x,\xi) = e^{2\pi \imath \Phi(x,\xi)}$ does not satisfy the
complementary low-rank property over the domain $X\times\Omega$ with
quadtree structures $T_X$ and $T_\Omega$. To fix this problem, the
polar butterfly algorithm introduces a scaled polar transformation on
$\Omega$:
\begin{equation}
  \xi=(\xi_1,\xi_2) = \frac{\sqrt{2}}{2}n p_1 \cdot (\cos{2\pi
    p_2},\sin{2\pi p_2}),
  \label{eq:polar}
\end{equation}
for $\xi\in \Omega$ and $p=(p_1,p_2)\in [0,1]^2$.  In the rest of this
section, we use $p$ to denote a point in the polar coordinate and $P$
for the set of all points $p$ transformed from $\xi\in\Omega$.  This
transformation gives rise to a new phase function $\Psi(x,p)$ in 
variables $x$ and $p$ satisfying
\begin{equation}
  \Psi(x,p) = \frac{1}{n}\Phi(x,\xi(p)) = \frac{\sqrt{2}}{2}\Phi\left(x, (\cos{2\pi p_2},\sin{2\pi p_2})\right) \cdot p_1,
\end{equation}
where the last equality comes from the fact that $\Phi(x,\xi)$ is
homogeneous of degree 1 in $\xi$. This new phase function $\Psi(x,p)$
is smooth in the entire domain $X\times P$ and the FIO in
\eqref{eq:fio} takes the new form
\begin{equation}
\label{eq:polarsum}
u(x) = \sum_{p\in P} e^{2\pi\imath n\Psi(x,p)}g(p),\quad x\in X.
\end{equation}

The transformation \eqref{eq:polar} ensures that $X\times P\subset
[0,1]^2\times [0,1]^2$. By partitioning $[0,1]^2$ recursively, we can
construct two quadtrees $T_X$ and $T_P$ of depth $L=\O(\log n)$ for
$X$ and $P$, respectively. The following theorem is a rephrased
version of Theorem 3.1 in \cite{fio09} that shows analytically the
complementary low-rank property of $e^{2\pi \imath n\Psi(x,p)}$ in the
$(X,P)$ domain.

\begin{theorem}
  \label{thm:pba}
  Suppose $A$ is a node in $T_X$ at level $\ell$ and $B$ is a node in
  $T_P$ at level $L-\ell$.  Given an {FIO} kernel function
  $e^{2\pi\imath n\Psi(x,p)}$ with a real-analytic phase function in
  the joint variables $x$ and $p$, there exist $\epsilon_0>0$ and
  $n_0>0$ such that for any positive $\epsilon\leq \epsilon_0$ and $n\geq n_0$,
  there exist $r_\epsilon$ pairs of functions $\{\alpha_t^{A,B}(x),
  \beta_t^{A,B}(p)\}_{1\leq t\leq r_\epsilon}$ satisfying that
  \begin{equation*}
    \left|e^{2\pi\imath n\Psi(x,p)} -
    \sum^{r_\epsilon}_{t=1}
    \alpha_t^{A,B}(x)\beta_t^{A,B}(p)
    \right|\leq \epsilon,
  \end{equation*}
  for $x\in A$ and $p\in B$ with $r_\epsilon\lesssim
  \log^4(1/\epsilon)$.
\end{theorem}

Based on Theorem \ref{thm:pba}, the polar butterfly algorithm
traverses upward in $T_\Omega$ and downward in $T_X$ simultaneously
and visits the low-rank submatrices $K_{A,B}=\{K(x_i,\xi_j)\}_{x_i\in
  A, \xi_j\in B}$ for pairs $(A,B)$ in $T_X\times T_P$. 
  The polar butterfly algorithm is
asymptotically very efficient: for a given input vector $g(p)$ for
$p\in P$, it evaluates \eqref{eq:polarsum} in $\O(N\log N)$ steps
using $\O(N)$ memory space. We refer the readers to \cite{fio09} for a 
detailed description of this algorithm.

\subsection{Factorization algorithm}
\label{sec:algopbf}


Combining the polar butterfly algorithm with the butterfly
factorization outlined in Section \ref{sec:gbf} gives rise to the
following polar butterfly factorization (PBF).
\begin{enumerate}
\item \emph{Preliminary.}
    Take the polar transformation of each point in $\Omega$
    and reformulate the problem
    \begin{equation}
    u(x) = \sum_{\xi\in \Omega} e^{2\pi\imath\Phi(x,\xi)}g(\xi),\quad x\in X,
    \end{equation}
    into
    \begin{equation}
    u(x) = \sum_{p\in P} e^{2\pi\imath n\Psi(x,p)}g(p),\quad x\in X.
    \end{equation}

\item \emph{Factorization.}
    Apply the two-dimensional {butterfly factorization}
    to the kernel $e^{2\pi\imath n\Psi(x,p)}$
    defined on a non-uniform point set in $X\times P$.
    The corresponding kernel matrix is approximated as
    \begin{equation}
         K\approx U^LG^{L-1}\cdots G^hM^h\left(H^h\right)^*\cdots
         \left(H^{L-1}\right)^*\left(V^L\right)^*.
    \end{equation}

\end{enumerate}

Since the polar butterfly factorization essentially applies the
original butterfly factorization to non-uniform point sets $X$ and
$P$, it has the same complexity as summarized in Table
\ref{tab:complexity}. Depending on the SVD procedure employed in the
middle level factorization, we refer to it either as PBF-m (when SVD
via random matrix-vector multiplication is used) or as PBF-s (when SVD
via random sampling is used).

\subsection{Numerical results}
\label{sec:numpbf}

This section presents two numerical examples to demonstrate the
efficiency of the polar butterfly factorization. The numerical results
were obtained in MATLAB on a server with 2.40 {GHz} {CPU} and 1.5 {TB}
of memory.


In this section, we denote by $\{u^p(x)\}_{x\in X}$ the results
obtained via the {PBF}. The relative error of the {PBF}
is estimated as follows, by comparing $u^p(x)$ with the exact values
$u(x)$.
\begin{equation}
  e^p = \sqrt{\cfrac{\sum_{x\in S}|u^p(x)-u(x)|^2}{\sum_{x\in
        S}|u(x)|^2}},
\end{equation}
where $S$ is a set of 256 randomly sampled points from $X$.


{\bf Example 1.}  The first example is a two-dimensional generalized
{Radon} transform that is an {FIO} defined as follows:
\begin{equation}
  \label{eq:example1fio}
  u(x) = \sum_{\xi\in\Omega} e^{2\pi \imath \Phi(x,\xi)}g(\xi),
\quad x\in X,
\end{equation}
with the phase function given by
\begin{equation}\label{eq:example1phi}
\begin{split}
\Phi(x,\xi) =& x\cdot \xi+\sqrt{c_1^2(x)\xi_1^2+c_2^2(x)\xi_2^2},\\
c_1(x) = & (2+\sin(2\pi x_1)\sin(2\pi x_2))/16,\\
c_2(x) = & (2+\cos(2\pi x_1)\cos(2\pi x_2))/16,
\end{split}
\end{equation}
where $X$ and $\Omega$ are defined in \eqref{eq:X} and
\eqref{eq:Omega}.  The computation in \eqref{eq:example1fio}
approximately integrates over spatially varying ellipses, for which
$c_1(x)$ and $c_2(x)$ are the axis lengths of the ellipse centered at
the point $x\in X$. The corresponding matrix form of
\eqref{eq:example1fio} is simply
\begin{equation}
  u=Kg, \quad K = (e^{2\pi\imath\Phi(x,\xi)})_{x\in X,\xi\in\Omega}.
  \label{eq:uKg}
\end{equation}

As $e^{2\pi \imath \Phi(x,\xi)}$ is known explicitly, we are able to
use the {PBF-s} (i.e., the one with random sampling in the middle
level factorization) to approximate the kernel matrix $K$ given by
$e^{2\pi \imath \Phi(x,\xi)}$. After the construction of the butterfly
factorization, the summation in \eqref{eq:example1fio} can be
evaluated efficiently by applying these sparse factors to $g(\xi)$.
Table \ref{tab:sec4example1} summarizes the results of this
example. 


\begin{table}[ht!]
\centering
\begin{tabular}{rcccc}
  \toprule
  $n,r$ & $\epsilon^p$ & $T_{f,p}(min)$ & $T_p(sec)$ & Speedup \\
\toprule
    64,6 & 2.46e-02 & 6.51e-01 & 2.37e-02 & 1.54e+02 \\
   128,6 & 7.55e-03 & 9.84e+00 & 2.30e-01 & 1.67e+02 \\
   256,6 & 5.10e-02 & 2.73e+01 & 6.23e-01 & 7.55e+02 \\
   512,6 & 1.46e-02 & 4.00e+02 & 7.88e+00 & 4.15e+02 \\
\toprule
    64,14 & 7.93e-04 & 7.34e-01 & 5.98e-02 & 8.72e+01 \\
   128,14 & 7.28e-04 & 1.17e+01 & 7.15e-01 & 4.28e+01 \\
   256,14 & 2.15e-03 & 3.93e+01 & 1.46e+00 & 2.86e+02 \\
   512,14 & 1.25e-03 & 5.63e+02 & 1.05e+01 & 3.35e+02 \\
\toprule
    64,22 & 6.96e-05 & 7.40e-01 & 8.24e-02 & 4.51e+01 \\
   128,22 & 7.23e-05 & 1.16e+01 & 1.04e+00 & 3.69e+01 \\
   256,22 & 2.44e-04 & 5.14e+01 & 5.94e+00 & 7.74e+01 \\
\bottomrule
\end{tabular}
\caption{Numerical results provided by the {PBF} with randomized sampling
    algorithm for the FIO in \eqref{eq:example1fio}.
    $n$ is the number of grid points in each dimension;
    $N=n^2$ is the size of the kernel matrix;
    $r$ is the max rank used in the low-rank approximation;
    $T_{f,p}$ is the factorization time of the {PBF};
    $T_d$ is the running time of the direct evaluation;
    $T_p$ is the application time of the {PBF}.
    The last column shows the speedup factor compared to the direct evaluation.
}
\label{tab:sec4example1}
\end{table}

{\bf Example 2.} The second example evaluates the composition of two
{FIOs} with the same phase function $\Phi(x,\xi)$. This is given
explicitly by
\begin{equation}\label{eq:example2mfio}
  u(x) =\sum_{\eta\in\Omega} e^{2\pi \imath \Phi(x,\eta)} \sum_{y\in
    X}e^{-2\pi \imath y\cdot \eta} \sum_{\xi\in\Omega} e^{2\pi \imath
    \Phi(y,\xi)}g(\xi), \quad x\in X,
\end{equation}
where the phase function is given in \eqref{eq:example1phi}.
The corresponding matrix representation is
\begin{equation}\label{eq:example2mfiomat}
u =KFKg,
\end{equation}
where $K$ is the matrix given in \eqref{eq:uKg} and $F$ is the matrix
representation of the discrete Fourier transform.  Under relatively
mild assumptions (see \cite{Theory} for details), the composition of
two FIOs is again an FIO. Hence, the kernel matrix
\begin{equation}
  \widetilde{K} := KFK
  \label{eq:KFK}
\end{equation}
of the product can be approximated by the butterfly factorization.
Notice that
the kernel function of $\widetilde{K}$ defined by \eqref{eq:KFK} is
not given explicitly. However, \eqref{eq:KFK} provides fast algorithms
for applying $\widetilde{K}$ and its adjoint through the fast
algorithms for $K$ and $F$. For example, the butterfly factorization
of Example 1 enables the efficient application of $K$ and $K^*$ in
$\O(N\log N)$ operations. Applying of $F$ and $F^*$ can be done by the
fast Fourier transform in $\O(N\log N)$ operations. Therefore, we can
apply the {PBF-m} (i.e., the one with random matrix-vector
multiplication) to factorize the kernel $\widetilde{K} = KFK$. Table
\ref{tab:sec4example2} summarizes the numerical results of this
example, the composing of two FIOs.

\begin{table}[ht!]
\centering
\begin{tabular}{rcccc}
  \toprule
  $n,r$ & $\epsilon^p$ & $T_{f,p}(min)$ & $T_p(sec)$ & Speedup \\
\toprule
    64,12 & 3.84e-02 & 6.22e+00 & 2.18e-02 & 3.34e+02 \\
   128,12 & 1.31e-02 & 3.86e+02 & 1.80e-01 & 4.25e+02 \\
\toprule
    64,20 & 2.24e-03 & 8.58e+00 & 3.04e-02 & 2.39e+02 \\
   128,20 & 2.23e-03 & 3.68e+02 & 3.60e-01 & 2.13e+02 \\
\bottomrule
\end{tabular}
\caption{Numerical results provided by the {PBF} with randomized {SVD}
    algorithm for the composition of FIOs
    given in \eqref{eq:example2mfiomat}.}
\label{tab:sec4example2}
\end{table}

{\bf Discussion.} The numerical results in Tables
\ref{tab:sec4example1} and \ref{tab:sec4example2} support the
asymptotic complexity analysis. When we fix $r$ and let $n$ grow, the
actually running time fluctuates around the asymptotic scaling since
the implementation of the algorithms differ slightly depending on
whether $L$ is odd or even. However, the overall trend matches well
with the $O(N^{3/2})$ construction cost and the $O(N\log N)$
application cost. For a fixed $n$, one can improve the accuracy by
increasing the truncation rank $r$. From the tables, one observes that
the relative error decreases by a factor of 10 when we increase the
rank $r$ by $8$ every time. In the second example, since the
composition of two {FIOs} typically has higher ranks compared to a
single FIO, the numerical rank $r$ used for the composition is larger
than that for a single FIO in order to maintain comparable accuracy.


\section{Multiscale Butterfly Factorization}
\label{sec:mbf}
In this section, we discuss yet another approach for constructing
butterfly factorization for the kernel $K(x,\xi)=e^{2\pi\imath
  \Phi(x,\xi)}$ with singularity at $\xi=0$. This is based on the
multiscale butterfly algorithm introduced in \cite{mba}.

\subsection{Multiscale butterfly algorithm}
\label{sec:mba}

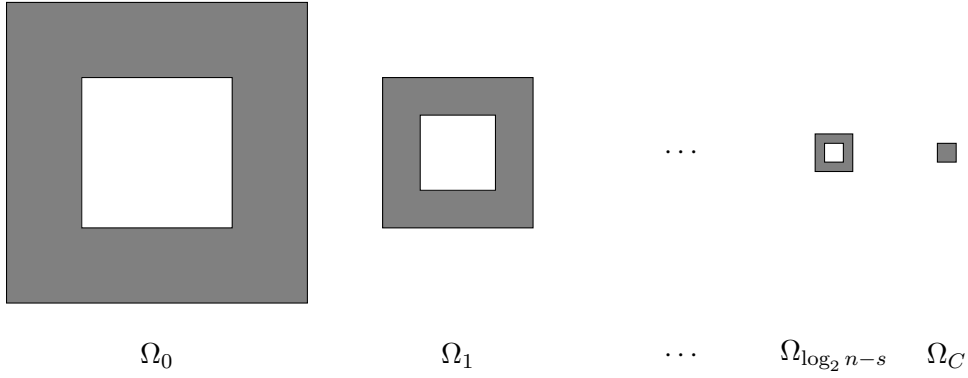
\begin{figure}[htb]
\centering
\begin{tikzpicture}[scale=1]
\coordinate (T) at (0,-2.7);

\fill[gray] (-2,-2) rectangle (2,2);
\draw[black] (-2,-2) rectangle (2,2);
\fill[white] (-1,-1) rectangle (1,1);
\draw[black] (-1,-1) rectangle (1,1);

\draw ($(T)$) node[rectangle] {$\Omega_0$};

\coordinate (A) at (4,0);
\fill[gray] ($(A)+(-1,-1)$) rectangle ($(A)+(1,1)$);
\draw[black] ($(A)+(-1,-1)$) rectangle ($(A)+(1,1)$);
\fill[white] ($(A)+(-.5,-.5)$) rectangle ($(A)+(.5,.5)$);
\draw[black] ($(A)+(-.5,-.5)$) rectangle ($(A)+(.5,.5)$);

\draw ($(A)+(T)$) node[rectangle] {$\Omega_1$};

\coordinate (A) at (7,0);
\draw ($(A)$) node[rectangle] {$\cdots$};
\draw ($(A)+(T)$) node[rectangle] {$\cdots$};

\coordinate (A) at (9,0);
\fill[gray] ($(A)+(-.25,-.25)$) rectangle ($(A)+(.25,.25)$);
\draw[black] ($(A)+(-.25,-.25)$) rectangle ($(A)+(.25,.25)$);
\fill[white] ($(A)+(-.125,-.125)$) rectangle ($(A)+(.125,.125)$);
\draw[black] ($(A)+(-.125,-.125)$) rectangle ($(A)+(.125,.125)$);

\draw ($(A)+(T)$) node[rectangle] {$\Omega_{\log_2 n -s}$};

\coordinate (A) at (10.5,0);
\fill[gray] ($(A)+(-.125,-.125)$) rectangle ($(A)+(.125,.125)$);
\draw[black] ($(A)+(-.125,-.125)$) rectangle ($(A)+(.125,.125)$);
\draw ($(A)+(T)$) node[rectangle] {$\Omega_C$};

\end{tikzpicture}
\caption{This figure shows the frequency domain decomposition of $\Omega$.
    Each subdomain $\Omega_t$, $t=0,1,\dots,\log_2 n-s$,
    is a corona subdomain
    and $\Omega_C$ is a small square subdomain covering the origin.}
\label{fig:domain-decomp}
\end{figure}

The key idea of the multiscale butterfly algorithm \cite{mba} is to
hierarchically partition the domain $\Omega$ into subdomains excluding
the singular point $\xi=0$.  This multiscale partition is illustrated
in Figure \ref{fig:domain-decomp} with
\begin{equation}\label{eq:domaindecomp}
  \Omega_t =
  \left\{(\xi_1,\xi_2):\frac{n}{2^{t+2}}<\max(|\xi_1|,|\xi_2|) \leq
  \frac{n}{2^{t+1}}\right\}\cap \Omega,
\end{equation}
for $t=0,1,\dots,\log_2 n-s$, $s$ is a small constant, and
$\Omega_C=\Omega\setminus\cup_t\Omega_t$. Equation
\eqref{eq:domaindecomp} is a corona decomposition of $\Omega$, where
each $\Omega_t$ is a corona subdomain and $\Omega_C$ is a square
subdomain at the center containing $\O(1)$ points.

The FIO kernel $e^{2\pi\imath\Phi(x,\xi)}$ satisfies the complementary
low-rank property when it is restricted in each subdomain $X\times
\Omega_t$. This observation is supported by the following theorem
rephrased from Theorem 3.1 in \cite{mba}. Here the notation
$\dist(B,0) = \min_{\xi\in B}\norm{\xi-0}$ is the distance between the
square $B$ and the origin $\xi=0$ in $\Omega$.

\begin{theorem}
\label{thm:mba}
Given an {FIO} kernel function $e^{2\pi\imath\Phi(x,\xi)}$ with a
real-analytic phase function $\Phi(x,\xi)$ for $x$ and $\xi$ away from
$\xi=0$, there exist a constant $n_0>0$ and a small constant
$\epsilon_0$ such that the following statement holds.  Let $A$ and $B$
be two squares in $X$ and $\Omega$ with sidelength $w_A$ and $w_B$,
respectively.  Suppose $w_Aw_B\leq 1$ and $\dist(B,0)\geq
\frac{n}{4}$.  For any positive $\epsilon\leq\epsilon_0$ and $n\geq
n_0$, there exist $r_\epsilon$ pairs of functions
$\{\alpha_t^{A,B}(x), \beta_t^{A,B}(p)\}_{1\leq t\leq r_\epsilon}$
satisfying that
\begin{equation*}
  \left|e^{2\pi\imath \Phi(x,\xi)} -
  \sum^{r_\epsilon}_{t=1}
  \alpha_t^{A,B}(x) \beta_t^{A,B}(\xi)
  \right|\leq \epsilon,
\end{equation*}
for $x\in A$ and $\xi\in B$ with $r_\epsilon\lesssim
\log^4(1/\epsilon)$.
\end{theorem}
According to the low-rank property in Theorem \ref{thm:mba}, the
multiscale butterfly algorithm rewrites \eqref{eq:fio} as a multiscale
summation,
\begin{equation}\label{eq:mba}
  u(x) = u_C(x) + \sum_{t=0}^{\log_2n-s} u_t(x) = \sum_{\xi\in
    \Omega_C}e^{2\pi\imath \Phi(x,\xi)}g(\xi)
  +\sum_{t=0}^{\log_2n-s}\sum_{\xi\in \Omega_t}e^{2\pi\imath
    \Phi(x,\xi)}g(\xi).
\end{equation}

For each $t$, the multiscale butterfly factorization
algorithm evaluates $u_t(x) = \sum_{\xi\in\Omega_t}
e^{2\pi\imath\Phi(x,\xi)} g(\xi)$ with a standard butterfly algorithm
such as the one that relies on the oscillatory {Lagrange}
interpolation on {Chebyshev} grid (see \cite{fio09}). The final piece
$u_C(x)$ is evaluated directly in $\O(N)$ operations. As a result, the
multiscale butterfly algorithm asymptotically takes $\O(N\log N)$
operations to evaluate \eqref{eq:mba} for a given input function
$g(\xi)$ for $\xi\in\Omega$. We refer the reader to \cite{mba} for the
detailed exposition.

\subsection{Factorization algorithm}
\label{sec:algombf}

Combining the multiscale butterfly algorithm with the butterfly
factorization outlined in Section \ref{sec:gbf} gives rise to the
following multiscale butterfly factorization (MBF).

\begin{enumerate}
\item \emph{Preliminary.}  Decompose domain $\Omega$ into subdomains
  as in \eqref{eq:domaindecomp}.  Reformulate the problem into a
  multiscale summation according to \eqref{eq:mba}:
  \begin{equation}
    \label{eq:msumd}
    K = K_C R_C + \sum_{t=0}^{\log_2n-s} K_t R_t.
  \end{equation}
  Here $K_C$ and $K_t$ are kernel matrices corresponding to $X\times
  \Omega_C$ and $X\times \Omega_t$. $R_C$ and $R_t$ are the
  restriction operators to the domains $\Omega_C$ and $\Omega_t$
  respectively.
\item \emph{Factorization.}  Recall that $L=\log_2 n$. For each
  $t=0,1,\dots,L-s$, apply the two-dimensional butterfly factorization
  on $K(x,\xi)=e^{2\pi \imath \Phi(x,\xi)}$ restricted in $X\times
  \Omega_t$.  Let $\widetilde{\Omega}_t$ be the smallest square that
  contains $\Omega_t$. Define $L_t = 2\lfloor (L-t)/2\rfloor$, where
  $\lfloor \cdot\rfloor$ is the largest integer less than or equal to
  a given number.  We construct two quadtrees $T_X$ and
  $T_{\widetilde{\Omega}_t}$ of depth $L_j$ with $X$ and
  $\widetilde{\Omega}_t$ being the roots, respectively. Applying the
  two-dimensional butterfly factorization using the quadtrees $T_X$
  and $T_{\widetilde{\Omega}_t}$ gives the $t$-th butterfly
  factorization:
  \begin{equation*}
    K_t \approx U_t^{L_t}G_t^{L_t-1}\cdots
    G_t^{\frac{L_t}{2}}
    M_t^{\frac{L_t}{2}}
    \left(H_t^{\frac{L_t}{2}}\right)^*\cdots
    \left(H_t^{L_t-1}\right)^*\left(V_t^{L_t}\right)^*.
  \end{equation*}
  Note that $1/4$ of the tree $T_{\widetilde{\Omega}_t}$ is empty and
  we can simply ignore the computation for these parts. This is a
  special case of non-uniform point sets. Once we have computed all
  butterfly factorizations, the multiscale summation in
  \eqref{eq:msumd} is approximated by
  \begin{equation}
    K \approx K_C R_C + \sum_{t=0}^{L-s}
    U_t^{L_t}G_t^{L_t-1}\cdots M_t^{\frac{L_t}{2}}\cdots
    \left(H_t^{L_t-1}\right)^*\left(V_t^{L_t}\right)^* R_t.
    \label{eq:f1}
  \end{equation}
\end{enumerate}

The idea of the hierarchical decomposition of $\Omega$ not only avoids
the singularity of $K(x,\xi)$ at $\xi=0$, but also maintains the
efficiency of the butterfly factorization.  The butterfly
factorization for the kernel matrix restricted in $X\times \Omega_t$
is a special case of non-uniform butterfly factorization in which the
center of $\Omega_t$ contains no point. Since the number of points in
$\Omega_t$ is decreasing exponentially in $t$, the operation and
memory complexity of the multiscale butterfly factorization is
dominated by the butterfly factorization of $K_t$ for $t=0$,
which is bounded by
the complexity summarized in Table \ref{tab:complexity}. Depending on
the SVD procedure in the middle level factorization, we refer this
factorization either as MBF-m (when SVD via random matrix-vector
multiplication is used) or as MBF-s (when SVD via random sampling is
used).

\subsection{Numerical results}
\label{sec:nummbf}

This section presents two numerical examples to demonstrate the
efficiency of the {MBF} as well.  The numerical results are obtained
in the same environment as the one used in Section \ref{sec:numpbf}.
Here we denote by $\{u^m(x),x\in X\}$ the results obtained via the
{MBF}. The relative error is estimated by
\begin{equation}
  e^m = \sqrt{\cfrac{\sum_{x\in S}|u^m(x)-u(x)|^2}{\sum_{x\in S}|u(x)|^2}},
\end{equation}
where $S$ is a set of 256 randomly sampled from $X$.  In the
multiscale decomposition of $\Omega$, we recursively divide $\Omega$
until the center part is of size 16 by 16.

{\bf Example 1. } We revisit the first example in Section \ref{sec:numpbf}
to illustrate the performance of the {MBF},
\begin{equation}\label{eq:sec5example1fio}
u(x) = \sum_{\xi\in\Omega} e^{2\pi \imath \Phi(x,\xi)}g(\xi),
\quad x\in X,
\end{equation}
with a kernel $\Phi(x,\xi)$ given by
\begin{equation}
\begin{split}
\Phi(x,\xi) =& x\cdot \xi+\sqrt{c_1^2(x)\xi_1^2+c_2^2(x)\xi_2^2},\\
c_1(x) = & (2+\sin(2\pi x_1)\sin(2\pi x_2))/16,\\
c_2(x) = & (2+\cos(2\pi x_1)\cos(2\pi x_2))/16,
\end{split}
\end{equation}
where $X$ and $\Omega$ are defined in \eqref{eq:X} and
\eqref{eq:Omega}.  Table \ref{tab:sec5example1} summarizes the results
of this example obtained by applying the {MBF-s}.

\begin{table}[ht!]
\centering
\begin{tabular}{rcccc}
  \toprule
  $n,r$ & $\epsilon^m$ & $T_{f,m}(min)$ & $T_m(sec)$ & Speedup \\
\toprule
    64,12 & 1.58e-02 & 4.48e-01 & 4.09e-02 & 1.13e+02 \\
   128,12 & 1.47e-02 & 5.64e+00 & 1.93e-01 & 2.02e+02 \\
   256,12 & 2.13e-02 & 2.16e+01 & 5.51e-01 & 9.26e+02 \\
   512,12 & 1.97e-02 & 2.97e+02 & 5.07e+00 & 6.45e+02 \\
\toprule
    64,20 & 5.51e-03 & 4.74e-01 & 6.11e-02 & 6.17e+01 \\
   128,20 & 4.27e-03 & 5.95e+00 & 5.01e-01 & 7.63e+01 \\
   256,20 & 1.68e-03 & 3.03e+01 & 2.51e+00 & 1.79e+02 \\
   512,20 & 2.02e-03 & 4.57e+02 & 1.14e+01 & 2.98e+02 \\
\toprule
    64,28 & 7.42e-05 & 7.18e-01 & 3.92e-02 & 6.23e+01 \\
   128,28 & 8.46e-05 & 1.23e+01 & 5.42e-01 & 7.43e+01 \\
   256,28 & 5.63e-04 & 6.73e+01 & 3.23e+00 & 1.43e+02 \\
   512,28 & 4.18e-04 & 7.20e+02 & 1.66e+01 & 2.14e+02 \\
\bottomrule
\end{tabular}
\caption{Numerical results provided by the {MBF} with the randomized
  sampling algorithm for the FIO given in \eqref{eq:sec5example1fio}.
  $n$ is the number of grid points in each dimension; $N=n^2$ is the
  size of the kernel matrix; $r$ is the max rank used in low-rank
  approximation; $T_{f,m}$ is the factorization time of the
  {MBF}; $T_d$ is the running time of the direct evaluation; $T_m$ is
  the application time of the {MBF}; $T_d/T_m$ is the speedup factor.}
\label{tab:sec5example1}
\end{table}

{\bf Example 2. } Here we revisit the second example in Section
\ref{sec:numpbf} to illustrate the performance of the {MBF}.  Recall
that the matrix representation of a composition of two {FIOs} is
\begin{equation}
  \label{eq:sec5example2mfiomat}
  u =\widetilde{K}g = KFKg,
\end{equation}
and that there are fast algorithms to apply $K$, $F$ and their
adjoints. Hence, we can apply the {MBF-m} (i.e., with the random
matrix-vector multiplication) to factorize $\widetilde{K}$ into the
form of \eqref{eq:f1}. Table \ref{tab:sec5example2} summarizes the
results.

\begin{table}[ht!]
\centering
\begin{tabular}{rcccc}
  \toprule
  $n,r$ & $\epsilon^m$ & $T_{f,m}(min)$ & $T_m(sec)$ & Speedup \\
  \toprule
  64,16 & 1.86e-02 & 4.05e+00 & 1.95e-02 & 4.23e+02 \\
  128,16 & 1.76e-02 & 1.27e+02 & 1.86e-01 & 4.17e+02 \\
  \toprule
  64,24 & 4.43e-03 & 5.37e+00 & 2.52e-02 & 3.27e+02 \\
  128,24 & 3.02e-03 & 1.79e+02 & 2.29e-01 & 3.40e+02 \\
  \bottomrule
\end{tabular}
\caption{MBF numerical results for the composition of FIOs
  given in \eqref{eq:sec5example2mfiomat}.}
\label{tab:sec5example2}
\end{table}

{\bf Discussion.} The results in Tables \ref{tab:sec5example1} and
\ref{tab:sec5example2} agree with the $\O(N^{3/2}\log N)$ complexity
analysis of the construction algorithm.  As we double the problem size
$n$, the factorization time increases by a factor 9 on average.  The
actual application time in these numerical examples matches the
theoretical operation complexity of $\O(N\log N)$.  In Table
\ref{tab:sec5example1}, the relative error decreases by a factor of 10
when the increment of the rank $r$ is 6.  In Table
\ref{tab:sec5example2}, the relative error decreases by a factor of 6
when the increment of the rank $r$ is 8.

\section{Conclusion}
\label{sec:conc}

We have introduced three multidimensional butterfly factorizations as
data-sparse representations of a class of kernel matrices coming from
multidimensional integral transforms. When the integral kernel
$K(x,\xi)$ satisfies the complementary low-rank property in the entire
domain, the butterfly factorization introduced in Section
\ref{sec:gbf} represents an $N\times N$ kernel matrix as a product of
$\O(\log N)$ sparse matrices. In the FIO case for which the kernel
$K(x,\xi)$ is singular at $\xi=0$, we propose two extensions: (1) the
polar butterfly factorization that incorporates a polar coordinate
transformation to remove the singularity and (2) the multiscale
butterfly factorization that relies on a hierarchical partitioning in
the $\Omega$ domain. For both extensions, the resulting butterfly
factorization takes $O(N\log N)$ storage space and $O(N \log N)$ steps
for computing matrix-vector multiplication as before.

The butterfly factorization for higher dimensions ($d>2$) can be
constructed in a similar way. For the {PBF}, one simply applies a
$d$-dimensional spherical transformation to the frequency domain
$\Omega$. For the {MBF}, one can again decompose the frequency domain
as a union of dyadic shells centered round the singularity at $\xi=0$.

{\bf Acknowledgments.}  This work was partially supported by the
National Science Foundation under award DMS-1328230 and the
U.S. Department of Energy's Advanced Scientific Computing Research
program under award DE-FC02-13ER26134/DE-SC0009409. H. Yang also
thanks the support from National Science Foundation under award 
ACI-1450372 and an AMS-Simons Travel Grant.

\bibliographystyle{abbrv} \bibliography{ref}

\begin{thebibliography}{10}

\bibitem{fio07}
E.~Cand{\`e}s, L.~Demanet, and L.~Ying.
\newblock Fast computation of {F}ourier integral operators.
\newblock {\em SIAM J. Sci. Comput.}, 29(6):2464--2493, 2007.

\bibitem{fio09}
E.~Cand{\`e}s, L.~Demanet, and L.~Ying.
\newblock A fast butterfly algorithm for the computation of {F}ourier integral
  operators.
\newblock {\em Multiscale Model. Simul.}, 7(4):1727--1750, 2009.

\bibitem{randsamp1}
B.~Engquist and L.~Ying.
\newblock A fast directional algorithm for high frequency acoustic scattering
  in two dimensions.
\newblock {\em Commun. Math. Sci.}, 7(2):327--345, 2009.

\bibitem{randsvd}
N.~Halko, P.~G. Martinsson, and J.~A. Tropp.
\newblock Finding structure with randomness: probabilistic algorithms for
  constructing approximate matrix decompositions.
\newblock {\em SIAM Rev.}, 53(2):217--288, 2011.

\bibitem{Theory}
L.~H\"{o}rmander.
\newblock {Fourier} integral operators. {I}.
\newblock {\em Acta Mathematica}, 127(1):79--183, 1971.

\bibitem{hu}
J.~Hu, S.~Fomel, L.~Demanet, and L.~Ying.
\newblock {A fast butterfly algorithm for generalized {Radon} transforms}.
\newblock {\em Geophysics}, 78(4):U41--U51, June 2013.

\bibitem{1dbf}
Y.~{Li}, H.~{Yang}, E.~{Martin}, K.~{Ho}, and L.~{Ying}.
\newblock {Butterfly Factorization}.
\newblock {\em Multiscale Model. Simul.}, 13(2):714--732, 2015.

\bibitem{mba}
Y.~{Li}, H.~{Yang}, and L.~{Ying}.
\newblock A multiscale butterfly algorithm for multidimensional {F}ourier
  integral operators.
\newblock {\em Multiscale Model. Simul.}, 13(2):614--631, 2015.

\bibitem{Rec2}
E.~Liberty, F.~Woolfe, P.-G. Martinsson, V.~Rokhlin, and M.~Tygert.
\newblock Randomized algorithms for the low-rank approximation of matrices.
\newblock {\em Proc. Natl. Acad. Sci. USA}, 104(51):20167--20172, 2007.

\bibitem{HMatrix}
L.~Lin, J.~Lu, and L.~Ying.
\newblock Fast construction of hierarchical matrix representation from
  matrix-vector multiplication.
\newblock {\em J. Comput. Phys.}, 230(10):4071--4087, 2011.

\bibitem{HSSMatrix}
P.~G. Martinsson.
\newblock A fast randomized algorithm for computing a hierarchically
  semiseparable representation of a matrix.
\newblock {\em SIAM J. Matrix Anal. Appl.}, 32(4):1251--1274, 2011.

\bibitem{mmd}
E.~Michielssen and A.~Boag.
\newblock A multilevel matrix decomposition algorithm for analyzing scattering
  from large structures.
\newblock {\em Antennas and Propagation, IEEE Transactions on},
  44(8):1086--1093, Aug 1996.

\bibitem{1dba}
M.~O'Neil, F.~Woolfe, and V.~Rokhlin.
\newblock An algorithm for the rapid evaluation of special function transforms.
\newblock {\em Appl. Comput. Harmon. Anal.}, 28(2):203--226, 2010.

\bibitem{fio13}
J.~Poulson, L.~Demanet, N.~Maxwell, and L.~Ying.
\newblock A parallel butterfly algorithm.
\newblock {\em SIAM J. Sci. Comput.}, 36(1):C49--C65, 2014.

\bibitem{wavemoth}
D.~S. Seljebotn.
\newblock Wavemoth-fast spherical harmonic transforms by butterfly matrix
  compression.
\newblock {\em The Astrophysical Journal Supplement Series}, 199(1):5, 2012.

\bibitem{sht}
M.~Tygert.
\newblock Fast algorithms for spherical harmonic expansions, {III}.
\newblock {\em J. Comput. Phys.}, 229(18):6181--6192, 2010.

\bibitem{Rec3}
F.~Woolfe, E.~Liberty, V.~Rokhlin, and M.~Tygert.
\newblock A fast randomized algorithm for the approximation of matrices.
\newblock {\em Applied and Computational Harmonic Analysis}, 25(3):335 -- 366,
  2008.

\bibitem{randsamp2}
H.~Yang and L.~Ying.
\newblock A fast algorithm for multilinear operators.
\newblock {\em Appl. Comput. Harmon. Anal.}, 33(1):148--158, 2012.

\bibitem{sft}
L.~Ying.
\newblock Sparse {F}ourier transform via butterfly algorithm.
\newblock {\em SIAM J. Sci. Comput.}, 31(3):1678--1694, 2009.

\end{thebibliography}

\end{document}